\documentclass[11pt]{amsart}
\usepackage{geometry}                % See geometry.pdf to learn the layout options. There are lots.
\usepackage{graphicx}
\usepackage{amssymb}
\usepackage{epstopdf}
\usepackage{tikz-cd}

\usepackage[T1]{fontenc}
\usepackage{svg}
\usepackage{relsize}
\usepackage{subcaption}

\usepackage{menukeys}
\usepackage{ wasysym }

\usepackage{enumerate}

\usepackage{verbatim}

\usepackage{bm}

\newcommand{\FF}{\mathbb F}

\newcommand{\II}{\mathbb I}

\newcommand{\cA}{\mathcal{A}}

\newcommand{\cC}{\mathcal{C}}
\newcommand{\cD}{\mathcal{D}}
\newcommand{\cE}{\mathcal{E}}
\newcommand{\cF}{\mathcal{F}}
\newcommand{\cG}{\mathcal{G}}
\newcommand{\cH}{\mathcal{H}}

\newcommand{\cO}{\mathcal{O}}

\newcommand{\cR}{\mathcal{R}}
\newcommand{\cS}{\mathcal{S}}
\newcommand{\cT}{\mathcal{T}}

\newcommand{\cV}{\mathcal{V}}
\newcommand{\cW}{\mathcal{W}}

\newcommand{\cZ}{\mathcal{Z}}

\newcommand{\inte}{\mbox{Int}}

\newcommand{\id}{\mathrm{Id}}

	\newcommand{\sfh}{SFH}
	
		\newcommand{\cfkr}{CFK_{\cR}}

	\newcommand{\bsd}{\widehat{\mathrm{BSD}}}
	
	\newcommand{\bsda}{\widehat{\mathrm{BSDA}}}
	\newcommand{\bsdd}{\widehat{\mathrm{BSDD}}}

		\newcommand{\cfdd}{\mathrm{CFDD}}

	\newcommand{\shi}{{SHI}}

	\newcommand{\shm}{\underline{{\mathrm{\bf SHM}}}}

	\newcommand{\bfd}{B\cF D}
	\newcommand{\invbfd}{BFD}
	
\newcommand{\bpsi}{\pmb{\psi}}

\newcommand{\bone}{\pmb{1}}

\newcommand{\bk}{{\bf k } }
\newcommand{\bj}{{\bf j} }

\newcommand{\bp}{{\bf p } }

\newcommand{\bx}{{\bf x } }
\newcommand{\by}{{\bf y } }

\newcommand{\bF}{{ \bf F} }
\newcommand{\bH}{{\bf H}}
\newcommand{\bM}{{ \bf M} }
\newcommand{\bT}{{\bf T}}
\newcommand{\bU}{{\bf U}}

\newcommand{\bcD}{\pmb{\cD}}

\newcommand{\bbox}{\boxtimes}
\newcommand{\cone}{{\text{Cone}}}

\theoremstyle{plain}
\newtheorem{thm}{Theorem}[section]
\newtheorem{thmquotes}{Moral}

\newtheorem{qn}[thm]{Question}
\newtheorem{conj}[thm]{Conjecture}

\newtheorem{lem}[thm]{Lemma}
\newtheorem{cor}[thm]{Corollary}
\newtheorem{prop}[thm]{Proposition}

\theoremstyle{definition}
\newtheorem{defn}[thm]{Definition}

\newtheorem{exmp}[thm]{Example}
\newtheorem{func}{Functor}

\theoremstyle{remark}
\newtheorem{rmk}[thm]{Remark}

\usepackage{scalerel,stackengine}
\stackMath
\newcommand\reallywidecheck[1]{%
\savestack{\tmpbox}{\stretchto{%
  \scaleto{%
    \scalerel*[\widthof{\ensuremath{#1}}]{\kern-.6pt\bigwedge\kern-.6pt}%
    {\rule[-\textheight/2]{1ex}{\textheight}}%WIDTH-LIMITED BIG WEDGE
  }{\textheight}% 
}{0.5ex}}%
\stackon[1pt]{#1}{\scalebox{-1}{\tmpbox}}%
}
\newcommand{\oppo}{\mathrm{op}}

\newcommand{\amod}{{}^{\cA}\mathsf{Mod}}
\newcommand{\tormod}{{}^{\cA(T^2)}\mathsf{Mod}}

\newcommand{\amoda}{{}^{\cA_1}\mathsf{Mod}^{\cA_2}}
\newcommand{\uamodau}{{}^{\cA_1}_{\hphantom{\cA}\mathfrak{u}}\mathsf{Mod}^{\cA_2}_{\mathfrak{u}}}

\newcommand{\amodada}{{}^{\cA_1}\mathsf{Mod}_{\cA_2}}

\newcommand{\amodt}{{}^{\cA}\mathsf{Mod}^{\cT}}

\newcommand{\Hamoda}{\mathsf{H}({ ^{\cA_1}\mathrm{Mod}^{\cA_2}})}
\newcommand{\Kamoda}{\mathsf{K}({ ^{\cA_1}\mathrm{Mod}^{\cA_2}})}

\newcommand{\contsuttor}{\mathrm{\mathsf{ContSutTor}}}
\newcommand{\contsut}{\mathrm{\mathsf{ContSut}}}
\newcommand{\diffsut}{\mathrm{\mathsf{DiffSut}}}
\newcommand{\diffsuttor}{\mathrm{\mathsf{DiffSutTor}}}

\newcommand{\diffborsut}{\mathsf{DiffBorSut}}
\newcommand{\diffborsutbtor}{\mathsf{DiffBorSut}\bm{\mathsf{Tor}}}

\newcommand{\diffborsuttor}{\mathrm{\mathsf{{DiffBorSutTor}}}}
\newcommand{\diffsutbtor}{\mathsf{DiffSut}\bm{\mathsf{Tor}}}

\newcommand{\contdiffsut}{\mathrm{\mathsf{{ContDiffSut}}}}
\newcommand{\contdiffsuttor}{\mathrm{\mathsf{{ContDiffSutTor}}}} 

\newcommand{\contdiffborsut}{\mathrm{\mathsf{{ContDiffBorSut}}}}

\newcommand{\ab}{\mathrm{\mathsf{Ab}}}
\newcommand{\moduvflip}{{\cR\mathsf{Mod}}^{\mathrm{flip}}}

\newcommand{\ob}{\mathrm{Ob}}

\newcommand{\psys}{\mathrm{\bf PSys}}

\newcommand{\csys}{\mathrm{\pmb \cC \bf Sys}}

\newcommand{\ama}{{}^{\cA_1} M^{\cA_2}}
\newcommand{\amt}{{ ^{\cA} M^{\cT}}}
\newcommand{\ana}{{ ^{\cA_1} {N}^{\cA_2}}}
\newcommand{\amada}{{ ^{\cA_1} M_{\cA_2}}}

\newcommand{\crC}{{\mathsf T}\cC}
\newcommand{\crpC}{{\mathsf T}^+\cC}
\newcommand{\TCh}{\mathsf{TCh}}
\newcommand{\TcC}{\crC}

\newcommand{\KC}{{\mathsf K} \cC}
\newcommand{\DC}{{\mathsf D} \cC}
\newcommand{\crCh}{{\mathsf{TCh}}}
\newcommand{\KCh}{{\mathsf{KCh}}}
\newcommand{\KcF}{\mathsf{K}\cF}
\newcommand{\Ch}{\mathsf{Ch}}

\newcommand{\dify}{\cmd}
\newcommand{\sunn}{{\normalfont \sun}}

\title{Bordered invariants of pairs of sutured manifolds with torus boundary}
\author{Tom Hockenhull}
\begin{document}
\maketitle
                                        % Activate to display a given date or no date
\begin{abstract}
We establish a framework for extending invariants of sutured manifolds to invariants of pairs of sutured manifolds who differ by attaching a basic slice along a torus boundary component. In the particular case of (bordered-)sutured Floer homology we show that this recovers the bordered-sutured Floer homology  of a corresponding bordered-sutured three-manifold --- though the framework is broad enough that it should be applicable in the setting of other invariants of sutured manifolds.\end{abstract}

To a category $\cS$ one can associate the corresponding \emph{arrow category} $\cS(2)$, whose objects are morphisms $f: A \rightarrow X$ and morphisms are commutative diagrams in $\cS$. There are variations on this theme: for instance, the \emph{category of pairs} of a concrete category takes as objects morphisms which are injective (as maps of sets). In more generality, to any category with two classes of morphism we can associate a corresponding category with objects given by one class of morphism, and morphisms given by commutative diagrams in $\cS$ with vertical arrows given by the other class.

A fairly natural question to ask is:
\begin{qn}What do invariants of and functors out of a category of pairs look like? How do they relate to functors and invariants of the underlying category?
\end{qn}

Of particular interest to us are categories where the objects are balanced sutured three-manifolds. To fix notation:
\begin{defn}A \emph{sutured three-manifold} is a pair $(M, \gamma)$, where
\begin{enumerate} \item $M$ is an oriented three-manifold-with-boundary, and
\item $\gamma \subset \partial M$ is a set of a set of pairwise disjoint annuli  whose core curves are given by oriented simple closed curves $s(\gamma) \subset \partial Y$.
\end{enumerate}
A sutured manifold $(M, \gamma)$ is a \emph{sutured submanifold} of a sutured manifold $(N, \delta)$ if $M \subset \inte(N)$.
\end{defn}
There are various choices for the morphism spaces, cf. \cite{bs:naturalitymonopole}. We are mostly interested in the category $\contdiffsut$ --- this has two classes of morphism:
\begin{enumerate} \item `gluing maps', i.e. whenever $(M, \gamma)$ is a sutured submanifold of $(N, \delta)$ there is a morphism $\psi: (M, \gamma) \rightarrow (N, \delta)$ corresponding to each isotopy class of contact structure upon $(N, \delta) - \inte(M, \gamma)$  which is convex on the boundary and whose dividing set induces the sutures $\gamma$ and $\delta$,  and
\item diffeomorphisms $\phi: (M, \gamma) \rightarrow (N, \delta)$.
\end{enumerate}

We are mostly interested in a category of pairs of the category $\contdiffsut$, which we call $\diffsut_2$, where the objects are given by gluing maps and morphisms by commutative diagrams where the vertical arrows are diffeomorphisms. In fact, we often further restrict our attention to a subcategory where the gluing map corresponds to \emph{attaching a basic slice}, which we call $\diffsuttor^\perp_2$.

It transpires that the category $\diffsuttor^\perp_2$ is equivalent to a category of bordered-sutured manifolds with a bordered torus boundary component, as defined in \cite{zarev:borderedsutured}. In particular, there is a well-known invariant $\bsd: \diffsuttor^\perp_2 \rightarrow \tormod$, where $\tormod$ is the category of type $D$ structures over an algebra associated to the torus $T^2$. This invariant, in a sense, `lifts' the invariant $\sfh: \diffsut \rightarrow \Ch$: there are functors $\pi_M$ and $\pi_N$ from $\diffsuttor^\perp_2$ to $\diffsut$, and $\iota_0 \bbox \cdot$ and $\iota_1 \bbox \cdot$ from $\tormod$ to $\Ch$, such that the diagrams
\[
\begin{tikzcd}
 f: M \rightarrow N\arrow[r,"\bsd"] \arrow[d,"\pi_M"]  & \bsd(f)\arrow[d,"\iota_0 \bbox \cdot"]  \\
M\arrow[r,"\sfh"]  & \sfh(M)
\end{tikzcd}\hspace{5em}
\begin{tikzcd}
 f: M \rightarrow N\arrow[r,"\bsd"] \arrow[d,"\pi_N"]  & \bsd(f)\arrow[d,"\iota_1 \bbox \cdot"]  \\
N\arrow[r,"\sfh"]  & \sfh(N)
\end{tikzcd}
\]
commute.

There are a number of other well-known invariants and functors out of $\diffsut$, cf. Subsection \ref{allthefunctors} below. The main philosophical take-away from this paper is:
\begin{thmquotes}\label{suturedtoborderedsutured}Most well-known functors $$\cF: \contdiffsut \rightarrow \Ch$$ should lift to functors $$\bfd: \diffsuttor_2^\perp \rightarrow \tormod.$$
\end{thmquotes}
We give explicit criteria for this to hold, and an explicit construction of the functor in this case.

It is often more difficult to construct a functor from $\contdiffsut$ to $\Ch$ than it is to construct an \emph{invariant} of $\contdiffsut$. Here, an invariant $F: \contdiffsut \rightarrow \cC$ is an assignment of an isomorphism class of object $[F(O)]$ in $\cC$ to each object $O$ of $\contdiffsut$. We can weaken our resolve suitably:
\begin{thmquotes}\label{suturedtoborderedsuturedinvariant}Most well-known invariants $$F: \contdiffsut \rightarrow \Ch$$ should lift to invariants $$\invbfd: \diffsuttor^\perp_2 \rightarrow \tormod.$$\end{thmquotes}

Besides this proselytizing and the accompanying theorems, the main result we show is:
\begin{thmquotes}\label{suturedfloertoborderedsuturedfloer}The bordered-sutured invariant $\bsd$ can be interpreted as an invariant of the category $\diffsuttor^\perp_2$ which lifts the invariant $\sfh$.
\end{thmquotes}

In the remainder of this introduction, we will make the statements of these morals more explicit, as Theorems.

\subsection{Pair-shaped}
We take a slightly different, but equivalent, tack to the category of pairs $\diffsut_2$.
\begin{defn}\label{suturedpair}A \emph{sutured pair} is a triple $((N, \delta), (A, \zeta), i)$ of a sutured manifold $(N, \delta)$, a contact three-manifold-with-boundary $(A, \zeta)$ with respect to which the boundary of $A$ is convex, and $i: (A, \zeta) \hookrightarrow (N, \delta)$ is an embedding which is surjective onto the boundary of $(N, \delta)$ --- such that the dividing set induced by $\zeta$ on the boundary of $N$ agrees with $\delta$.

We often suppress the inclusion map $i$ from notation.
\end{defn}

A sutured pair is equivalent to a pair in the category $\contsut$ as above, cf. also in \cite{bs:naturalitymonopole} --- the complement $N - A$ is naturally a sutured submanifold of $(N, \delta)$ with sutures induced by the contact structure $\zeta$.

\begin{defn}The \emph{diffeomorphism category of sutured pairs} is the category $\diffsut_2$ whose objects are sutured pairs as above, and whose morphisms between $((N, \delta), (A, \zeta), i)$ and $((N', \delta'), (A', \zeta'), i')$ are pairs $(\psi_N, \psi_A)$ of diffeomorphisms $\psi_N: (N, \delta) \rightarrow (N', \delta')$ and $\psi_A: (A, \zeta) \rightarrow (A', \zeta')$ such that $\psi_N \circ i = i' \circ \psi_A$ and $\psi_A^* \zeta' = \zeta$.
\end{defn}

We often restrict our attention to the subcategory of $\diffsut_2$ which we call the \emph{category of orthogonal pairs $\diffsuttor_2^{\perp}$}. Here, we require $(N, \delta)$ to have a torus boundary component $T$ and $(A, \zeta)$ to be a \emph{basic slice}, cf. \cite[Section 2.1]{honda:tightcontacti} --- that is, a minimally twisting contact structure upon $T^2 \times I$ which induces the correct sutures on $T \subset \partial N$. In order for this to be a bona fide pair we must extend $(A, \zeta)$ to surject onto the boundary of $N$; we do this by taking $(A, \zeta) = (\partial N \times I, \zeta)$ and choosing $\zeta$ to restrict to a product contact structure on $(\partial N - T) \times I$ and a basic slice on $T \times I$.

Topologically, in this case, the manifold $(M, \gamma) = (N, \delta) - (A, \zeta)$ is the same as $(N, \delta)$, and away from the distinguished torus boundary component $T$ we have that $\gamma|_{\partial Y -T} = \delta|_{\partial Y -T}$. On the component $T$ the sutures are parallel copies of a curve which intersects each component of $\delta|_{T}$ transversally in a single point: in this case we say that $\delta$ and $\gamma$ are \emph{orthogonal} or, indeed, that $(M, \gamma)$ and $(N,\delta)$ are \emph{orthogonal}.

\begin{figure}
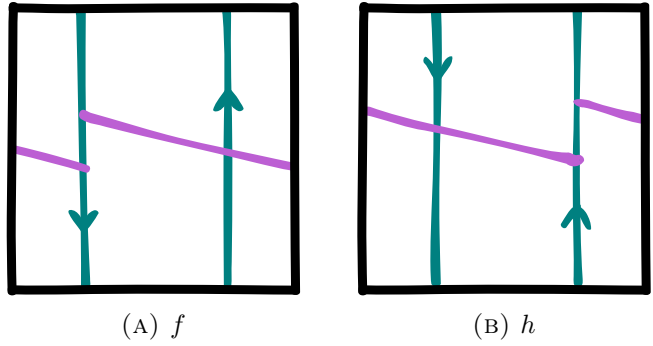

  \begin{minipage}{0.3\linewidth}\centering
    \scalebox{1}{\includesvg[pretex=\relscale{1}]{bypass1}}    \subcaption{ $f$}
  \end{minipage}
    \begin{minipage}{0.3\linewidth}\centering
   \scalebox{1}{\includesvg[pretex=\relscale{1}]{bypass2}}    \subcaption{ $h$}
  \end{minipage}

  \caption{The two choices of bypass arc.\label{fig:bypasses}}%

\end{figure}

There are precisely two basic slices which induce any fixed pair of orthogonal sutures, and these can be seen in terms of a \emph{positive} and \emph{negative} bypass on a thickened torus $T^2 \times I$ endowed with an $I$-invariant contact structure. In essence this means taking the sutures upon $T^2 \times \{1\}$ and attaching half an overtwisted disc along an arc which begins and ends upon one of the sutures, and intersects the other in a single point. There are two choices, $\zeta_f$ and $\zeta_h$ corresponding with arcs $f$ and $h$ as shown in Figure \ref{fig:bypasses}, cf. \cite[Figure 6]{stipsiczvertesi:oninvariants}.
 
We will say that the two pairs $((N, \delta), (A, \zeta_f))$ and $((N, \delta), (A, \zeta_h))$ are \emph{conjugate}.

In essence, then, an orthogonal pair is determined by a sutured manifold $(N, \delta)$ with preferred torus boundary component $T$ together with a choice of bypass arc $f \subset T$. We write $((N, \delta), f)$ for this pair, and, having chosen a fixed $f$, we will write $((N, \delta), h)$ for the conjugate pair. We'll also write $f(N, \delta)$ for the sutured manifold $(M, \gamma) = (N, \delta) - (A, \zeta)$.

\subsubsection{Triangles}
The category $\contdiffsut$ is supposed to be triangulated, with details in unpublished work by Honda. To avoid straying too far from the course, we settle here for describing a set of `preferred triangles'.

\begin{defn}A \emph{preferred triangle} in $\contdiffsut$ is a triple $((M_i, \gamma_i), f_i)_{i=0}^{i=3}$ of sutured manifolds $M_i$ and gluing maps $f_i: M_i \rightarrow M_{i+1}$, such that the $(M_i, \gamma_i)$ differ locally as shown in Figure \ref{fig:triangleofmaps}, and the maps $f_i$ correspond with attaching a bypass along the arcs shown. (Compare, say, \cite[Theorem 1.21]{bs:khovanovtrefoils}.)
\begin{figure}
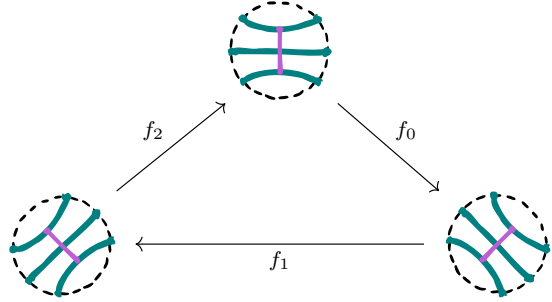

\centering
   \[ 
\begin{tikzcd}
 &\scalebox{1}{\includesvg[pretex=\relscale{1}]{bypasstriangle}}    \arrow[
      start anchor={[shift={(0.4cm,0.1cm)}]}, % fine tune start of arrow
      end anchor={[shift={(0.35cm,-0.15cm)}]} % fine tune end of arrow
      ]{dr}{f_0} \\
{\raisebox{-1cm}{\rotatebox{45}{{\includesvg[pretex=\relscale{1}]{bypasstriangle}}}}}   \arrow[
      start anchor={[shift={(-0.4cm,-0.15cm)}]}, % fine tune start of arrow
      end anchor={[shift={(-0.35cm,0.1cm)}]} % fine tune end of arrow
      ]{ur}{f_2}&& \rotatebox{-45}{\scalebox{1}{\includesvg[pretex=\relscale{1}]{bypasstriangle}}}  \arrow[
      start anchor={[shift={(0.15cm,-0.12cm)}]}, % fine tune start of arrow
      end anchor={[shift={(-0.15cm,-0.12cm)}]} % fine tune end of arrow
      ]{ll}{f_1} \end{tikzcd}
\] 
   
  \caption{The local picture in a preferred triangle \label{fig:triangleofmaps}}%

\end{figure}

\end{defn}

In particular, given any pair $((N, \delta), f)\in \ob(\contdiffsut_2^\perp)$ there is a corresponding triangle $\Delta_f$ whose corresponding sutured manifolds $(M_i, \gamma_i)$ satisfy that $(M_0, \gamma_0)=(N, \delta)$ and subsequent sutured manifolds are obtained by attaching basic slices of the same sign as $f$.

There is, of course, another \emph{conjugate triangle} corresponding to the conjugate pair $((N, \delta), h)$: here, the corresponding sutured manifolds are the same, but the slices attached are always of the opposite sign to $f$.

\begin{rmk}One obstacle to declaring $\contdiffsut$ triangulated is the absence of obvious shift functors; addressing this is in and of itself an interesting question --- cf. \cite{huang:bypassplanefields}, \cite{me:iteratedcones}.
\end{rmk}

\subsubsection{Bordered-sutured manifolds with bordered torus boundary component}
We are also interested by another class of three-manifold, with a priori different boundary data.
\begin{defn}
A \emph{balanced sutured manifold with bordered torus boundary component} is a tuple $(Y, \Gamma, \phi)$, where:
\begin{enumerate} \item $Y$ is an oriented three-manifold-with-boundary and distinguished torus boundary component $T$,
\item $\gamma \subset \partial Y - T$ is a set of pairwise disjoint annuli $\gamma \subset \partial Y - T$ whose core curves are given by oriented simple closed curves $s(\gamma) \subset \partial Y - T$, and
\item $\phi: F(\cZ) \rightarrow T$ is an orientation-preserving homeomorphism, from the pointed surface $F(\cZ)$ associated to the pointed matched circle $\cZ = \cZ(T^2)$ --- see Figure \ref{fig:arcsurface}, and \cite[Figure 11.2]{LOT}

\begin{figure}
  \begin{minipage}{0.3\linewidth}\centering
    \scalebox{1}{\includesvg[pretex=\relscale{1}]{puncturedtorus}}  \subcaption{\label{fig:arcsurface} The surface $F(\cZ(T^2))$.}
  \end{minipage}
    \begin{minipage}{0.3\linewidth}\centering
   \scalebox{1}{\includesvg[pretex=\relscale{1}]{torustobypass}}   \subcaption{\label{fig:bypassdata} Bypass data.}
   \end{minipage}
 \label{fig:arcdiagrams}
\end{figure}

\end{enumerate}
and we require that every component of $(\partial Y - T) - \inte(\gamma)$ is oriented coherently with respect to the set $s(\gamma)$. Writing $R_+(\gamma)$ for the components of this whose orientation agrees with that induced by the orientation of $M$, and $R_-(\gamma)$ for the complement in $(\partial Y - T) $ of $R_+(\gamma)$, we say a sutured manifold with bordered torus boundary component is \emph{balanced} if $\chi(R_+(\gamma)) = \chi(R_-(\gamma))$.

A diffeomorphism of balanced sutured manifolds with bordered torus boundary component $(Y, \gamma, \phi)$ and $(Y', \gamma', \phi')$ is a diffeomorphism $\psi: Y \rightarrow Y'$ such that  $\phi' = \psi \circ \phi$.

Together these form a category, $\diffsutbtor$.\end{defn}

We've already been dealing with these in disguise:
\begin{prop}There is an equivalence of categories $\sunn: \diffsuttor^\perp_2 \rightarrow \diffsutbtor$.
\end{prop}
\begin{proof}Fix a balanced sutured manifold with bordered torus boundary component $(Y, \gamma, \phi)$, and consider Figure \ref{fig:bypassdata}. We let $\gamma_{T} \subset F(\cZ)$ be the pair of sutures pictured, and let $f \subset F(\cZ)$ by the arc shown. The homeomorphism $\phi$ realises these as a balanced sutured manifold $(M, \gamma \cup \gamma_{T})$ with a preferred torus boundary component $T$ together with a bypass arc $f \subset T$ --- specifying an orthogonal pair as in the previous section. To a diffeomorphism $\psi: (Y, \Gamma, \phi) \rightarrow (Y', \Gamma', \phi')$ we associate the corresponding diffeomorphism from $((M, \gamma \cup \gamma_{T^2}, f)$ to $((M', \gamma' \cup \gamma'_{T^2}, f')$.

This functor, $\ast$,  is clearly full, faithful and essentially surjective --- so gives an equivalence of categories. We let $\sunn$ denote the inverse functor.
\end{proof}

\subsection{Functors and invariants of $\contdiffsut$, statement of theorems}\label{allthefunctors}
In this paper, we tend to assume that we are only given a suitable \emph{weak} invariant of $\contdiffsut$ rather than a bona fide functor out of $\contdiffsut$. Nevertheless, it is more straightforward (linguistically, at least) to understand the lay of the land if we work with a functor  $\cF: \contdiffsut \rightarrow \Ch$ (or $\ab$).

In this setting, if we want the functor $\cF$ to respect triangles in the sense that any favoured triangle of $\contsut$ is sent to an exact triangle, we have to instead work with a functor to $\KCh$ to even make sense of triangles being exact. This is easily remedied: given such a functor there is an associated functor $\mathsf{K}\cF: \contdiffsut \rightarrow \KCh$ obtained by passing to the homotopy category.

In this setting, one usually uses the \emph{triangle detection lemma}, cf. \cite{os:bdc,seidel:fukaya}. The purpose of this lemma is to build, for any preferred triangle $\Delta = ((M_i, \gamma_i), f_i)_{i=0}^{i=2}$ of $\contdiffsut$ a pair of quasi-inverse chain maps
\begin{align*}\Phi^{f_0}_\Delta: \cone(\cF(f_0)) &\rightarrow \cF(M_2, \gamma_2)\\
\Phi^{f_0}_\nabla: \cF(M_2, \gamma_2) &\rightarrow \cone(\cF(f_0)).\end{align*}

This then ensures that the corresponding triangle $\mathsf{K}\cF(\Delta)$ is exact, when working over a field, say, so that the homotopy and derived categories of chain complexes agree. A choice of such maps for every preferred triangle will be called a set of \emph{biwitnesses} for functor $\cF$ --- as they `witness' the fact that $\mathsf{K}\cF$ is triangulated. A pair $(\cF, \cW)$ of functor and biwitnesses is a \emph{witnessed functor}.

The important moral here is that `in nature', it is reasonable to expect that for any triangulated functor of the form $\mathsf{K}\cF$ we have a canonical choice of witnesses for $\cF$: to show that $\mathsf{K}\cF$ is triangulated in the first place we essentially have to build \emph{some} set of witnesses for $\cF$ --- the only leap of faith is in expecting these to be suitably independent of choices.

\subsubsection{Carriage returns and typewriters}
The category $\contdiffsut$ has a curious and essential feature related to its `triangulated' structure. As we have seen, given a sutured pair $((M, \gamma), f) \in \ob(\diffsuttor_2^\perp)$ there are two associated triangles $\Delta_f$ and $\Delta_h$ who share the \emph{same} sutured manifolds.

It follows immediately that for any functor $\cF: \contdiffsut \rightarrow \Ch$ such that the corresponding $\KcF: \contdiffsut \rightarrow \KCh$ respects preferred triangles, the two complexes $\cone(\cF(f))$ and $\cone(\cF(h))$ must be chain homotopy equivalent. Moreover, if we choose a set of biwitnesses for $\cF$, this gives an explicit choice of chain homotopy equivalence
$$\Phi^{h}_\nabla \circ \Phi^{f}_\Delta: \cone(\cF(f)) \rightarrow \cone(\cF(h)),$$
for any orthogonal pair $((M; \gamma), f) \in \diffsuttor^\perp_2$. 
(It also gives an explicit choice of chain homotopy inverse to these by swapping the roles of $f$ and $h$.)

We call such a choice of chain homotopy equivalences a set of \emph{carriage returns} for $\cF$. (A choice of the homotopy inverses is called a set of \emph{carriage departures}, and a choice of both is called a set of \emph{carriages} for $\cF$.) In fact, we extract this set-up, and define a dg category of \emph{typewriters of chain complexes} $\TCh$ whose objects are tuples
$$\bM = (M_0, M_1; D_f, D_h; D_{CR}),$$
where $M_0, M_1$ are chain complexes, $D_f, D_h: M_0 \rightarrow M_1$ are chain maps, and where $D_{CR}: \cone(D_f) \rightarrow \cone(D_h)$ is a chain map called a \emph{carriage return}.

We prove:
\begin{thm}\label{liftingupfunctors}To any biwitnessed functor $(\cF, \cW)$ valued in $\Ch$ there is an associated functor $\cF_{2}^\cW: \diffsuttor^\perp_2 \rightarrow \TCh$ which lifts $\cF$.
\end{thm}

As aforementioned, we work mostly with \emph{weak invariants} as opposed to functors. With this in mind, we show:
\begin{thm}\label{carriagesforsfhexist}There is a well-defined set of carriage returns for the weak invariant $\sfh: \contdiffsut \rightarrow \Ch$, so that $\sfh$ lifts to an invariant $B\sfh D: \diffsuttor_2^\perp \rightarrow \TCh$.
\end{thm}

This seems somewhat esoteric in its utility --- the punchline, however, is:
\begin{thm}There is an equivalence of categories $\dify: \TCh \rightarrow \tormod$, where $\cA(T^2)$ is the algebra associated to a torus by \cite{LOT}.\label{typewritersaretyped} 
\end{thm}

There is nothing new under the sun:
\begin{thm}The following diagram commutes:
$$\begin{tikzcd}\diffsuttor_2^\perp \arrow[r, "B\sfh D"] \arrow[d, swap, "{\normalfont \sunn}"]& \TCh \arrow[d, "\dify"]\\
\diffsutbtor \arrow[r, "CFD"] & \tormod.
\end{tikzcd}$$
\end{thm}

\subsubsection{Extendability}
The attentive reader will note that we didn't require a carriage return to be a chain homotopy equivalence, as is suggested by the information witnesses give us. Requiring the existence of a homotopy inverse $D_{CD}$ (which we call a \emph{carriage departure}) is a stronger condition which we call being \emph{semi-extendable}. A \emph{semi-extended} typewriter is a typewriter together with a specific choice of carriage departure $D_{CD}$ which is a homotopy inverse to $D_{CD}$. These form a subcategory $\mathsf{T}^+\mathsf{Ch}$ of $\TCh$, and it is clear from our discussion that a biwitnessed functor naturally yields a semi-extended typewriter.

There is a corresponding notion of semi-extendability for type $D$ structures over $\cA(T^2)$, and we can refine Theorem \ref{typewritersaretyped}:
\begin{thm}The equivalence of categories $\dify: \TCh \rightarrow \amod$ sends $\mathsf{T}^+\mathsf{Ch}$ to the subcategory of semi-extendable type $D$ structures.
\end{thm}

\subsubsection{Weak invariants}
As we've discussed, it is often considerably harder to define a functor from $\contdiffsut$ than it is to simply define a \emph{weak invariant} of $\contdiffsut$, valued in $\cC$, i.e. an assignment of an isomorphism class of object in $\cC$ to every object of $\contdiffsut$. The usual strategy for this is to make a number of choices of \emph{auxiliary data} associated to an object of $\contdiffsut$, and define an associated object of $\cC$ using this: one then verifies that any two different choices of data yield isomorphic objects in $\cC$. To pass to a functor, one must verify that there is a \emph{canonical} isomorphism between the objects associated to any two choices of data.

With this in mind, most of our work concerns weak invariants, largely in order that one can still obtain a \emph{weak} invariant of $\diffsuttor^\perp_2$ without having to embark on proving full naturality results. The essential idea is to make a series of compatible choices of auxiliary data associated to objects in $\cC$, morphisms between them, and exact triangles ---  to define an object of $\TcC$ ---  then understand and define homotopy equivalences between objects corresponding with different choices of auxiliary data for the same object of $\diffsuttor^\perp_2$.

\subsection{Functors in the wild}
This theory is all well and good, but, of course, one would like to \emph{apply} it. Sadly, the current literature does not quite contain the requisite ingredients to apply Theorem \ref{liftingupfunctors}. Showing the functoriality of an invariant is difficult in the first place, and, in particular, the naturality of witnesses or carriage returns has not been addressed anywhere --- nor have the requisite invariance maps for pairs.

Nevertheless, we are close. We posit that for a triangulated category $\cC$ any functor $\cF: \contdiffsut \rightarrow \cC$ which has been shown to respect preferred triangles by the methods of the exact triangle detection lemma should have corresponding natural witness maps.

\begin{func}\label{suturedfunctor}The functor $\sfh(-)$ is an example of a functor from $\contsuttor$ to the category of chain complexes over $\FF_2$. Here, the morphisms for diffeomorphisms are defined and shown to be functorial in \cite{jtz:naturality}, and gluing morphisms are defined in \cite{hkm:gluing} --- and re-evaluated in \cite{jz:contacthandles} where they are shown to be natural.

A proof that the associated functor to the homotopy category respects triangles is somewhat non-standard: it is shown in \cite{zevv} using techniques of bordered-sutured Floer homology which avoid the need for the use of the exact triangle detection lemma, and by \cite{ls:equivalencegluing} the maps involved are equivalent to those defined in \cite{hkm:gluing}.

In particular, we do not automatically obtain candidate witness maps. Nevertheless, we show in Section \ref{borderedsuturedsection} that there is a good choice choice of carriage returns for the associated weak invariant, such that Theorem \ref{carriagesforsfhexist} holds.
\end{func}

\begin{func}\label{suturedmonopolefunctor}In \cite{km:knotssuturesexcision}, Kronheimer and Mrowka define a monopole invariant for sutured manifolds, $\shm(M)$. Baldwin and Sivek show in \cite{bs:naturalitymonopole} that this gives a functor --- valued in a category they call \emph{projectively transitive systems}, $\psys$. With this in mind, one hopes for a corresponding \emph{projectively transitive system of type $D$ structures} to be obtained by the methods in this paper.

In \cite[Theorem 5.2]{bs:contactsuturedmonopole} it is further shown that there are corresponding bypass morphisms, and that these satisfy bypass exact triangles. The proof of the latter boils down constructing closures of the sutured manifolds in question, identifying the bypass maps with cobordisms between these closures, and applying the \emph{surgery exact triangle} in monopole Floer homology, shown in \cite[Section 5]{kmos:monopoleslensspacesurgeries}; this uses the exact triangle detection lemma, hence should give a candidate set of witness morphisms dependent upon these choices, and therefore a type $D$ structure which we conjecture is an invariant of the pair.
\end{func}

\begin{func}\label{suturedinstantonfunctor}In \cite{km:knotssuturesexcision}, an invariant for sutured manifolds $\shi(M)$ is defined. Bypass morphisms are constructed in \cite{bs:instantoncontact, bs:khovanovtrefoils}; moreover, in \cite[Theorem 1.21]{bs:khovanovtrefoils} these are shown to respect triangles, along similar lines to the functor $\shm(M)$. The proof here (as with $\shm(M)$) constructs auxiliary closures of the sutured manifolds in question, identifying the maps in the triangle with cobordism maps and applying the surgery exact triangle in instanton Floer homology, cf. \cite{scaduto:instantonskhovanov} which uses the triangle detection lemma --- thus giving us a candidate set of witness morphisms, and therefore a type $D$ structure associated to these choices which we conjecture is an invariant of the pair.
\end{func}

\subsection{Future directions}
Besides the obvious intended purpose of the theory constructed in this paper --- to enable one to define bordered versions of other invariants --- there are a number of other interesting avenues one could follow.

\subsubsection{Computing with typewriters}Whilst somewhat abstract, the theory in this paper gives a good vantage point from which to compute bordered-sutured invariants $\bsd$ for three-manifolds with torus boundary components. 

The strength here is that the invariant $\sfh$ can be interpreted as a variant of the `hat' version of \emph{knot Floer homlogy}, for knots in sutured three-manifolds. Much work has gone into understanding versions of this for closed three-manifolds: in particular, there are formulae allowing one to compute the corresponding invariants $\sfh(N, \delta)$ and $\sfh(M, \gamma)$ for a pair, when $(N, \delta)$ represents a nullhomologous knot $K$ in a three-manifold $Y$, in terms of invariants $\cfkr(Y,K)$, cf. \cite{heddenlevine:surgery, hedden:whitehead} --- these are generally called `mapping cone theorems' in the theory.

In the sequel \cite{me:iteratedcones}, we use similar techniques as inspiration to give a formula for the bordered invariants $\bsd$ of three-manifolds with a torus boundary component, in terms of analogues of the invariant $\cfkr$. This uses in an essential manner the theory built in this paper: the setting of typewriters and carriage returns is an invaluable framework in which to do these calculations.

As an example of this sort of set-up, let $\cR = \FF[U,V]/({UV}=0)$. The category of \emph{differential $\cR$-modules with a flip} has objects $(M, \varphi)$ where $M$ is a differential $\cR$-module, and $\varphi: \cone(U) \rightarrow \cone(V)$ is a chain homotopy equivalence --- it has morphisms given by a pair $(f, H)$ of a $U, V$ equivariant homomorphism $f$ of differential $\cR$-modules $f$, and a chain homotopy $H$ between the composites $\phi' \circ f$ and $f \circ \phi$.

It is straightforward to see:
\begin{prop}There is a functor $\div: \moduvflip \rightarrow \crCh$.
\end{prop}

The invariant $\cfkr$ is a functor $\mathsf{Knot}_* \rightarrow \moduvflip$ from the category of pointed knots $\mathsf{Knot}_*$ with morphisms given by diffeomorphisms, cf. \cite{zemke:linkcobordisms}. It follows immediately that:
\begin{prop}There is a functor $\bsd_\infty: \mathsf{Knot}_* \rightarrow \tormod$.
\end{prop}
\begin{proof}This is given by the composite $\bsd_\infty : = \dify \circ \div \circ \cfkr$.
\end{proof}

\subsubsection{Naturality}The bordered-sutured invariants $\bsd$ are not presently known to give a functor, just a weak invariant; the constructions outlined in this paper give explicit criteria for a bordered invariant $BFD$ defined according to the framework below to give a strong invariant --- or functor. It would be interesting to use this structure to establish functoriality properties of the bordered invariants $\bsd$ --- given the work of \cite{jtz:naturality}, half of the requisite conditions to be checked are already known to be true.

\subsubsection{Different categories}
In this paper, we work with the category $\diffsuttor_2^\perp$, where we are preoccupied by torus boundary components. It is natural to look for a generalisation to this where we work with a more general category, like the category of pairs $\diffsut_2$ --- i.e. we consider corresponding bordered invariants whose boundary has any genus.

A different generalisation would be to consider a wider class of vertical morphism in the commutative diagrams defining $\diffsuttor_2^\perp$: for instance, we could consider $\contdiffsuttor_2^\perp$ whose vertical morphisms are allowed, themselves, to be given by attaching a basic slice or bypass. Perhaps more interestingly, one could consider the broadest category $\mathsf{SutTor}_2^\perp$ where we also allow vertical morphisms to be given by cobordisms of sutured manifolds which are fixed near the boundary --- this would allow one to define cobordism maps in the bordered-sutured theory.

\subsubsection{2-categories}
There is quite a striking and mysterious interpretation of a typewriter $\bM \in \ob(\TCh)$: namely, as a \emph{cell of a 2-category}. We refer the reader to \cite{maclane:categories} for an overview of the theory of 2-categories.

We let $\cC$ be the 2-category where: \begin{enumerate}\item objects are chain complexes over $\FF_2$,
\item 1-morphisms are chain maps, and \item 2-morphisms $\alpha: f \Rightarrow g: M \rightarrow N$ are chain maps $\alpha: \cone(f) \rightarrow \cone(g)$.\end{enumerate}
In order to actually specify a 2-category, we need to define horizontal and vertical composition of 2-cells. Vertical composition is straightforward: the composite of $\alpha: f \Rightarrow g$ and $\alpha': g \Rightarrow h$ is just $\alpha' \circ \alpha$. Horizontal composition is considerably more elaborate, and we defer this to a later stage where we have more notation, i.e. Section \ref{typewritersaretwocategory}. There appears to be no mention of this 2-category in the literature: it generalises the usual 2-category of chain complexes.

In particular, a 2-morphism in $\cC$ gives the same data as a typewriter in $\ob(\TCh)$, i.e. a type $D$ structure over $\cA(T^2)$. With this 2-category structure in mind, a natural conjecture is:
\begin{conj}\label{twocategoryconjecture}There is a 2-category structure upon $\contdiffsut$ such that $\diffsuttor_2^\perp$ is a 2-subcategory of this, and such that for pairs of conjugate basic slices there is a two-morphism between the corresponding gluing maps.
\end{conj}
If this were the case, we would obtain:
\begin{cor}Any functor of 2-categories  $\cF: \contdiffsut \rightarrow \cC$ assigns a typewriter to any pair in $\diffsuttor_2^\perp$, and therefore a type $D$ structure over $\cA(T^2)$ to any balanced bordered-sutured three-manifold with bordered torus boundary component.
\end{cor}

\subsection{Acknowledgements}
I'd like to thank Robert Lipshitz, Jake Rasmussen and Steven Sivek for useful and interesting conversations. I'm particularly grateful to Marco Golla for many useful conversations and his patience, and to Marco Marengon for useful conversations --- and to both for comments on a draft of this paper.

Much of the research in this paper was carried out at the Max Planck Institute for Mathematics.

\section{Background on type $D$ structures and bimodules}
We give a whirlwind tour of some of the algebra concerning type $D$ structures and bimodules that we will need later. Mostly this is for the fixing of notation: the only (slightly) novel section is regarding \emph{semi-extendable $DD$ bimodules}.

Throughout, as elsewhere in this paper, we suppress gradings.

\subsection{$DD$ bimodules}
Let $\cA = (A, \partial_A)$ be a dg algebra, and $\bk$ its ring of idempotents.
\begin{defn}[{\cite[Definition 2.2.55]{lot:bimodules}}] A \emph{(left) type $D$ structure over $\cA$} is a $\bk$-module $M$, together with a `structure' map
$\delta^1: M \rightarrow A \otimes M$
satisfying that $(\II \otimes \delta^1) \circ \delta^1 = 0$. Often one writes ${ ^\cA N}$ for the data of the pair $(M, \delta^1)$.

If $\cA = \cA_1 \otimes_\cR \cA_2^{\oppo}$ then we say that $(M, \delta^1)$ is an \emph{DD}-bimodule over $(\cA_1, \cA_2)$; here, we think of the structure map as $\delta^1:M \rightarrow A_1 \otimes M \otimes A_2$ and often write $\ama$ for the pair $(M, \delta^1)$.

There is a notion of a \emph{morphism} $\phi^1: \ama \rightarrow \ana$ of type $D$ structures and type $DD$ bimodules; these consist of a map $\phi^1: M_1 \rightarrow A_1 \otimes M_2 \otimes A_2$ satisfying a structure relation. As such, there is a dg category $\uamodau$ whose objects are $DD$ bimodules and morphisms are such morphisms.

One adjective to describe $DD$ bimodules is to be \emph{bounded}; that is, $\ama$ is bounded if for every $x \in M$ there is some $n$ such that $\delta^i(x) = 0$ for all $i >n$. We use this to define well-behaved subcategories of type $DD$ bimodules: if $\cA_1$ and $\cA_2$ be are algebras then denote by $\amoda$ the full subcategory of $\uamodau$ which are homotopy equivalent to a bounded type $DD$ structure.

As $\amoda$ is a dg category, one can form the associated \emph{derived category} $\Hamoda$ by passing to homology of morphism spaces. There is also an appropriate notion of a \emph{homotopy} between maps of type $D$ structures, taking the form of a map $\nu: M_1 \rightarrow A_1 \otimes M_2 \otimes A_2$ satisfying certain conditions; this yields an appropriate notion of two type $DD$ bimodules being \emph{homotopy equivalent} if there is a homotopy invertible morphism between them, and, in much the same way as in the category of chain complexes $\Ch$, allows us to form the \emph{homotopy category} $\Kamoda$ whose objects are those of $\amoda$ and morphisms are homotopy equivalence classes of morphisms.
\end{defn}

It's worth noting here that this nomenclature is controversial: oftentimes what we call the derived category is also called the homotopy category in dg category literature. Nevertheless, under our boundedness assumptions on $\amoda$ these two constructions agree here, cf. \cite[Corollary 2.4.4]{lot:bimodules} --- allowing us to neatly sidestep any contentiousness. We work almost exclusively in the homotopy category in this paper.

\begin{exmp}A particular case of this construction is where $\cA_1 = \cA_2$ is the trivial (unital) algebra; in this case, we recover the category of chain complexes over $\FF_2$, for which we write $\Ch$ with associated derived category $\mathsf{H}\Ch$ and homotopy category $\KCh$.\end{exmp}

\subsection{Cones}
An awful lot of this paper will be spent talking in cones.
\begin{defn}Let $\ama$ and $\ana$ be type $DD$ bimodules, and $\phi: \ama \rightarrow \ana$ be a morphism of type $DD$ bimodules. There is a $DD$ bimodule defined by $${}^{\cA_1}\cone(\phi^1)^{\cA_2} := M \oplus N,$$ equipped with the structure map $$\delta_{\cone}^1(x, y) := (\delta^1(x), \delta^1(y) + \phi^1(x)).$$
\end{defn}
This enables us to endow $\Kamoda$ with the structure of a triangulated category:
\begin{defn}A triangle of $DD$ bimodules $({ ^{\cA_1}M_i^{\cA_2}}, \phi_i)_{i = 0}^2$ in $\Kamoda$ is \emph{exact} if it is isomorphic to a triangle of the form

\[
\begin{tikzcd}
 & { ^{\cA_1}M_0^{\cA_2}}\arrow[dr,"\phi_0"] \\
{ ^{\cA_1}\cone(\phi_0)^{\cA_2}} \arrow[ur,"\phi_2"] &&\arrow[ll,"\phi_1"] {}^{\cA_1}M_1^{\cA_2}
\end{tikzcd}.
\]

\end{defn}

\begin{rmk}This is probably the worst consequence of our omission of gradings; to strictly have triangulated structure we have to keep track of gradings and introduce shift functors. The troubled reader is invited to replace `triangulated category' with `category with preferred triangles' and interpret the above triangle accordingly.
\end{rmk}

\subsection{$A_\infty$ modules and box tensor products}
The counterpart to type $D$ structures as defined above are \emph{$\cA_\infty$ modules}, cf. eg. \cite{LOT}. A hybrid of the two are \emph{$DA$ bimodules}.
\begin{defn}Let $\cA_1, \cA_2$ be dg algebras over rings $\bk$ and $\bj$, respectively. A \emph{type $DA$ bimodule} over $\cA_1$ and $\cA_2$ is a $(\bk, \bj)$-bimodule $N$ together with maps
$$\delta^1_{1+j}: N \otimes B \rightarrow A^{\otimes j} \otimes N,$$
satisfying a compatibility condition, cf. \cite[Definition 2.2.43]{lot:bimodules}. We will often write $\amada$ to denote a $DA$-bimodule over $\cA_1$ and $\cA_2$.

A \emph{morphism} of $DA$ bimodules $f^1: {}^{\cA_1}M_{\cA_2} \rightarrow {}^{\cA_1}N_{\cA_2}$ is a collection of maps 
$f^1_{1+j}: M \otimes A_2 \rightarrow A_1 \otimes N$
satisfying a compatibility condition; there are appropriate notions of homotopy of morphisms, homotopy equivalence and boundedness such that we can define a category $\amodada$ whose objects are those $DA$ bimodules homotopy equivalent to a bounded $DA$ bimodule, with morphisms given as above.
\end{defn}

When $\cA_1$ is trivial, this recovers the notion of a $\cA_\infty$ module over $\cA_2$. There is an appropriate notion of a tensor product of a $DA$ bimodule with a $DD$ bimodule:
\begin{defn}Let $\amada$ and ${ ^{\cA_2}N^{\cA_3}}$ be a $DA$ bimodule over $\cA_1, \cA_2$ and a $DD$ bimodule over $\cA_2, \cA_3$ respectively. Then there is an associated $DD$ bimodule over $\cA_1, \cA_3$, written as $\amada \bbox { ^{\cA_2}N^{\cA_3}}$.
\end{defn}

One can also tensor morphisms:
\begin{defn}Let $\phi: { ^{\cA_2}N_1^{\cA_3}} \rightarrow { ^{\cA_2}N_1^{\cA_3}}$ be a morphism of $DD$ bimodules. Then there is a corresponding morphism
$$\II \bbox \phi:  \amada \bbox { ^{\cA_2}N_1^{\cA_3}} \rightarrow  \amada \bbox { ^{\cA_2}N_1^{\cA_3}},$$
cf. \cite[Figure 5]{lot:bimodules}.
\end{defn}

This tells us that $\amada \bbox \cdot: { ^{\cA_2}\mathsf{Mod}^{\cA_3}} \rightarrow { ^{\cA_1} \mathsf{Mod}^{\cA_3}} $ constitutes a functor.

\subsection{Coefficient maps and generalised coefficient maps}

We will work a lot with a specific algebra over the course of this paper:
\begin{defn}
Let $\cT$ be the path algebra of the quiver
$$ 
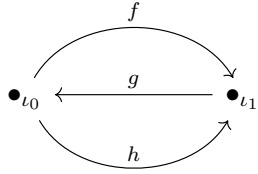
\begin{tikzcd}[column sep=60pt]
\bullet_{\iota_0} \arrow[bend left=60]{r}{f} \arrow[bend right=60]{r}{h}& \bullet_{\iota_1}  \arrow[swap]{l}{g}
\end{tikzcd}$$
after taking the quotient by the relations $gf = hg = 0$.
\end{defn}
This is isomorphic to the algebra $\cA(T^2)$ associated to the pointed torus by \cite{LOT}. For any type $D$ structure over the algebra $\cT$ there is a fairly established convention on how to describe the map $\delta^1$: in terms of \emph{coefficient maps}, cf. \cite{LOT}. We will need to extend this convention to a more general $DD$ bimodule over a pair $(\cA, \cT)$.

\begin{defn}Let $\bk$ denote the ring of idempotents of an algebra $\cA \otimes \cT$, and $M$ be $\bk$-module. A \emph{set of coefficient maps} is the data of maps $D_\emptyset^0, D_\emptyset^1, D_f, D_g, D_h, D_{fg}, D_{gh},$ and $D_{fgh}$ from $M$ to $\cA \otimes M$ such that the expressions
\begin{alignat*}{3}
D_\emptyset \circ D_f  + D_f \circ D_\emptyset &  \qquad D_\emptyset \circ D_g  + D_g \circ D_\emptyset  &\qquad  D_\emptyset \circ D_h  + D_h \circ D_\emptyset \end{alignat*}
\begin{alignat*}{2}
  D_{fg}\circ D_\emptyset + D_\emptyset \circ D_{fg} + D_g \circ D_f  &\qquad  D_{gh} \circ D_\emptyset + D_\emptyset \circ D_{gh} + D_h \circ D_g  
\end{alignat*}
all evaluate to zero.
\end{defn}

\begin{prop}\label{coefftodd}Let $(M, D_*)$ be a module and set of coefficient maps. The formula
$$\delta^1 = (D_\emptyset^0, D_\emptyset^1) \otimes \bone + \sum_{t \in \cT'} D_{t} \otimes t$$
gives $M$ the structure of a $DD$ bimodule over $(\cA, \cT)$.

Conversely, any $DD$ bimodule over $(\cA, \cT)$ uniquely determines a set of coefficient maps.
\end{prop}
\begin{proof}This is a matter of unwinding the definitions.
\end{proof}

\begin{rmk}Viewing a $DD$ bimodule in this manner is essentially the same as viewing it as a pre-twisted complex valued in the category of type $D$ structures over $\cA$, rather than the category of vector spaces, cf. \cite[Section 31]{seidel:fukaya}.\end{rmk}

There is a fairly formal extension of this construction, cf. \cite{LOT}.
\begin{defn}\label{semiextendable}Let $\amt$ be a $DD$ bimodule over the pair $(\cA, \cT)$. A set of \emph{generalised coefficient maps} compatible with $\amt$ is a set of maps $D_I$ for any proper cyclic interval $I \subset \{0, 1, 2, 3\}$ and further maps $D_{0123}, D_{1230}, D_{2301}, D_{3012}$. These are required to satisfy
$$\sum_{I = J \cup K | J < K} D_K \circ D_J = 0$$
for $<$ the cyclic ordering on $\{0,1,2,3\}$ and $I$ any proper cyclic interval in $\{0,1,2,3\}$. There are four more compatibility conditions required to be satisfied:
\begin{align*}D \circ D_{0123} + D_3 \circ D_{012} + D_{23} \circ D_{01} + D_{123} \circ D_0 + D_{0123} \circ D &= \II\\
D \circ D_{1230} + D_0 \circ D_{123} + D_{30} \circ D_{12} + D_{230} \circ D_1 + D_{1230} \circ D &= \II\\
D \circ D_{2301} + D_1 \circ D_{230} + D_{01} \circ D_{23} + D_{301} \circ D_2 + D_{2301} \circ D &= \II\\
D \circ D_{3012} + D_2 \circ D_{301} + D_{12} \circ D_{30} + D_{012} \circ D_3 + D_{3012} \circ D &= \II.
\end{align*}
We will say that a $DD$ bimodule $\amt$ over $(\cA, \cT)$ is \emph{semi-extendable} if there exist a set of generalised coefficient maps compatible with $\amt$, in the sense that there exists a set of generalised coefficient maps with $D_1 = D_f$, $D_{2} = D_g$, $D_3 = D_h$, $D_{12} = D_{fg}$, $D_{23} = D_{gh}$ and $D_{123} = D_{fgh}$.
\end{defn}

\begin{rmk}The notion of a $DD$ bimodule being semi-extendable is an approximation to that of (in the case of $\cA$ the trivial algebra) being \emph{extendable} in the sense of \cite{HRW}. In a little more detail: for a type $D$ structure over $\cT$ to be \emph{extendable} is to have a corresponding a curved type $D$ structure over an algebra $\widetilde{\cA}$ which contains algebra elements $\rho_I$ for any cyclic interval of $\{0,1,2,3\}$ --- the corresponding coefficient maps (in the appropriate sense) give maps $D_I$ which satisfy the compatibility conditions above.
\end{rmk}

\subsection{Bordered-Sutured manifolds, pairs and pairing}
We discuss the generalisation from categories of sutured manifolds to categories of bordered-sutured manifolds, and recall useful properties of bordered-sutured manifolds which we will need in this paper.
\subsubsection{Categories of Bordered-Sutured manifolds}
A bordered-sutured manifold $(Y, \Gamma, \cZ, \phi)$  is a three-manifold-with-boundary $Y$ together with some parameterisation and partitioning of its boundary, cf. \cite[Section 3]{zarev:borderedsutured}; a \emph{diffeomorphism} of bordered-sutured manifolds is a diffeomorphism of underlying three-manifolds-with-boundary which respects the parameterisation of the boundary. Together, these give a category $\diffborsut$.

A special case of this is where the portion of the data $\cZ$ corresponding to a boundary component $\partial_iY$ of $Y$ is empty, then we will say that component is \emph{balanced sutured}. In this manner a bordered-sutured manifold generalises a balanced sutured manifold: if $\cZ$ is empty then $(Y, \Gamma)$ is a balanced sutured manifold. In the other extreme, in the case where a portion of the data $\cZ$ corresponding to a boundary component $\partial_i Y$ is a pointed matched circle in the sense of \cite{LOT} we shall say that component is \emph{bordered}.

Generalising the definition of $\contdiffsut$ to bordered-sutured manifolds is somewhat fiddly in the most generality; cf. \cite{ls:equivalencegluing}. The following suffices for us:
\begin{defn}The category $\contdiffborsut$ has as objects bordered-sutured three-manifolds; there are two types of morphism:
\begin{enumerate} \item `sutured gluing maps', i.e. the set of sutured gluing maps between $(Y, \Gamma, \cZ, \phi)$ and $(Y', \Gamma', \cZ', \phi')$ is nonempty unless $Y$ embeds into $Y'$ such that the complement $A = Y' - \inte(Y)$ has boundary coinciding with the disjoint union of the balanced sutured boundary components of $Y$ and $Y'$. In this case there is a morphism for every isotopy class of contact structure upon $A$ for which the contact structure is convex on the boundary and induces the corresponding sutures for $Y$ and $Y'$ as dividing set.
\item diffeomorphisms $\phi: (M, \gamma) \rightarrow (N, \delta)$.
\end{enumerate}
\end{defn}

The restrictions on the sutured gluing maps are perhaps a little more reasonable when viewed from the perspective of pairs:
\begin{defn}
\label{borderedsuturedpair}A \emph{bordered-sutured pair} is a triple $((Y, \Gamma, \cZ, \phi), (A, \zeta), i)$ of a bordered-sutured manifold $(Y, \Gamma, \cZ, \phi)$, a three-manifold-with-boundary $A$ and isotopy class of contact structure $\zeta$ with respect to which the boundary of $A$ is convex, and $i: A \hookrightarrow Y$ is an embedding which is surjective onto the balanced sutured boundary of $(Y, \Gamma, \cZ, \phi)$ and whose image is disjoint from all non balanced sutured boundary components of $Y$ ---  and such that the dividing set induced by $\zeta$ on the balanced-sutured boundary of $Y$ agrees with the sutures in $\Gamma$.

We often suppress the inclusion map $i$ from notation.

The \emph{diffeomorphism category of bordered-sutured pairs} is the category $\diffborsut_2$ whose objects are bordered-sutured pairs as above, and whose morphisms between two triples  $((Y, \Gamma, \cZ, \phi), (A, \zeta), i)$ and $((Y', \Gamma', \cZ', \phi'), (A', \zeta'), i')$  are given by pairs of diffeomorphisms $\psi_Y: (Y, \Gamma, \cZ, \phi) \rightarrow(Y', \Gamma', \cZ', \phi')$ and $\psi_A: (A, \zeta) \rightarrow (A', \zeta')$ such that $\psi_Y \circ i = i' \circ \psi_A$ and $\psi_A^* \zeta' = \zeta$.
\end{defn}

With this in mind, the definition of the category $\diffsuttor_2^\perp$ generalises naturally to the subcategory $\diffborsuttor_2^\perp$ of $\diffborsut_2$ --- i.e. we again require $(A, \zeta)$ to be a (product extension of a) basic slice. Oncemore, in this context, a pair is specified by a pair $(Y, f)$ of a bordered-sutured manifold $Y$ with a balanced sutured torus boundary component $T$, and  a bypass arc $f \subset T$ --- we also write $f(Y)$ for the bordered-sutured manifold obtained by attaching a basic slice corresponding with $f$, and write $(Y, h)$ for the corresponding conjugate pair. A set of preferred triangles in $\contdiffborsut$ can also be defined analogously to those in $\contdiffsut$ --- in particular, to every bordered-sutured pair $(Y, f)$ we may associate a corresponding preferred triangle $\Delta_f$, with conjugate preferred triangle $\Delta_h$.

The definition of the category $\diffsutbtor$ readily generalises to this setting too --- giving a category $\diffborsutbtor$ whose objects are bordered-sutured manifolds with a specified bordered boundary component which is parameterised by the pointed matched circle associated to a torus. The definition of the functor $\sunn$ extends likewise, giving an equivalence of categories $\sunn: \diffborsuttor_2^\perp \rightarrow \diffborsutbtor$ .

\subsubsection{Bordered-Sutured diagrams and invariants}
To the parameterisation and partitioning of the boundary of a bordered-sutured manifold $(Y, \Gamma, \cZ, \phi)$ one associates a differential graded algebra $\cA(\cZ)$ --- if the boundary data splits as a disjoint union $\cZ_1 \cup \cZ_2$ along collections of boundary components then this algebra splits as a tensor product $\cA(\cZ_1) \otimes \cZ(\cZ_2)$; the algebra is trivial when the boundary component is balanced sutured.

A bordered-sutured manifold is determined by a \emph{bordered-sutured Heegaard diagram}. In \cite{zarev:borderedsutured}, it is shown how to associate to a bordered-sutured Heegaard diagram $\cH$ a cohort of invariants: of particular interest to us are the invariants $\bsd(\cH)$ and $\bsda(\cH)$. In the former case this is a type $D$ structure over the algebra $\cA(\cZ)$: as we have seen, where $(Y, \Gamma, \cZ, \phi)$ has multiple boundary components $\partial_1 Y \cup \partial_2 Y$ this takes the form of a type $D$ structure over $\cA(\cZ_1) \otimes \cA(\cZ_2)$, i.e. it is a $DD$ bimodule $\bsdd(\cH)$ over the algebra $(\cA(\cZ_1), \cA(\cZ_2))$. In the same circumstance, the invariant $\bsda(\cH)$ is a $DA$ bimodule over the algebra $(\cA(\cZ_1), \cA(\cZ_2))$.

These are shown to be a weak invariant of the underlying bordered-sutured manifold, in the sense of Section \ref{weakinvariants}, i.e. that the isomorphism classes of $\bsda(\cH)$ and $\bsdd(\cH)$ depend only on $(Y, \Gamma, \cZ, \phi)$ rather than the specific choice of $\cH$. More specifically, any two choices $\cH$ and $\cH'$ can be connected by some sequence of \emph{isotopies, handleslides and stabilisations}, cf. \cite[Theorem 7.8]{zarev:borderedsutured} and \cite[Section 6.3]{LOT}, which induce homotopy equivalences between $\bsdd(\cH)$ and $\bsdd(\cH')$ (resp. $\bsda(\cH)$ and $\bsda(\cH')$). In this circumstance one often writes $\bsda(Y)$ and $\bsdd(Y)$ for these isomorphism classes.

The invariants are also shown to obey a \emph{pairing theorem}: namely, $\cH = \cH_1 \cup_{\cZ_2} \cH_2$ is formed by gluing together two bordered-sutured Heegaard diagrams then their associated algebraic objects glue accordingly:
$$\bsdd(\cH) \cong \bsda(\cH_1) \bbox \bsdd(\cH_2).$$

\section{Weak and strong invariants}
In this section, we explore a general framework for defining invariants and functors from categories which rely on making choices as an intermediate step. Essentially, we expand the situation in \cite{jtz:naturality} to encompass that in \cite{bs:naturalitymonopole} --- also taking into account morphisms, and preferred triangles.

\subsection{Auxiliary data, weak and strong invariants}
The general schematic we use for building an invariant of a category $\cS$ or a functor out of $\cS$ is to first associate auxiliary data which encodes the objects, morphisms and any other data we are interested in, and then to this associate objects, morphisms, etc. of another category $\cC$ --- then check that these are suitably well-defined regardless of choice.
\subsubsection{Auxiliary data}
We make more explicit what we mean by \emph{auxiliary data} in this context.

\begin{defn}Let $\cS$ be a category, and $O \in \ob(\cS)$. A \emph{set of auxiliary data for $O$} is a connected graph $\cG(O)$ with vertex set $\cV(M)$ and edge set $\cE(O)$. %We call a vertex $\cD \in \cV(O)$ an \emph{auxiliary datum for $O$} and an edge $e \in \cE(M)$ an \emph{invariance map}.

A \emph{set of auxiliary object data for $\cS$} is a choice of auxiliary data $\cG(O)$ for every $O \in \ob(\cS)$.
\end{defn}
% Of particular interest to us in this paper is the setting where $\cS$ is the category $\bsutt$, and each object is of the form $O = (Y, \gamma)$.

We can extend the definition to take morphisms into account:
\begin{defn}Let $\cS$ be a category, and $\cG$ be a set of auxiliary object data for $\cS$. If $O, O'$ are objects of $\cS$ and $f: O \rightarrow O'$ is a morphism, then a \emph{set of auxiliary data for $f$} is a connected nonempty subgraph $\cG(f)$ of $\cG(O) \times \cG(O')$ --- with vertices $\cV(f) \subset \cV(O) \times \cV(O')$ and edges $\cE(f) \subset \cE(O) \times \cE(O')$.

A \emph{set of auxiliary data} for the category $\cS$ is a set of auxiliary object data for $\cS$ together with a choice of set of auxiliary data $\cG(f)$ for any morphism of $\cS$.
\end{defn}

Finally, if $\cS$ is a category with preferred triangles, we'd like to bear this structure in mind:
\begin{defn}Let $\cS$ be a category with preferred triangles, and $\cG$ be a set of auxiliary data for $\cC$. For any preferred triangle $\Delta = (O_i, f_i)_{i=0}^{i=2}$ in $\cS$, a \emph{set of auxiliary data for $\Delta$} is a connected nonempty subgraph $\cG(\Delta)$ of $\cG(O_0) \times \cG(O_1) \times \cG(O_2)$ such that the projection of $\cG(\Delta)$ to any cyclically consecutive pair $\cG(O_i) \times \cG(O_{i+1})$ is a subgraph of $\cG(f_i)$.

A \emph{set of auxiliary data for $\cS$ with triangles} is a set of auxiliary data for $\cS$ together with a choice of auxiliary data $\cG(\Delta)$ for any preferred triangle $\Delta$ of $\cS$. Usually we drop the \emph{with triangles} descriptor and assume that a set of auxiliary data for any category with preferred triangles is one with triangles.
\end{defn}

\subsubsection{Weak invariants}\label{weakinvariants}
With this formalism in hand, we can recover a generalisation of the notion of a \emph{weak Heegaard invariant} \cite[Definition 2.24]{jtz:naturality}.

\begin{defn}Let $\cS$ be a category, and $\cG = \cG(\cS)$ be a set of auxiliary object data for $\cS$. For any other category $\cC$, a \emph{weak $\cG$ invariant of $\ob(\cS)$ (valued in $\cC$)} is a morphism of graphs $F: \cG(\cS) \rightarrow \cC$ such that for any edge $e \in \cE(\cS)$, the image of $e$ is an isomorphism in $\cC$.
\end{defn}

As in loc. cit., if $F: \cG(\cS) \rightarrow \cC$ is a weak $\cG$ invariant then to each object $O \in \ob(\cS)$ we can associate an isomorphism class of object of $\cC$ --- i.e. this isomorphism class is an invariant of $O$.

We will need a further generalisation of this definition, to a more general set of auxiliary data:
\begin{defn}\label{weakinvariant}Let $\cS$ be a category, and $\cG = \cG(\cS)$ be a set of auxiliary data for $\cS$. For any other category $\cC$, a \emph{weak $\cG$ invariant of $\cS$ (valued in $\cC$)} is a weak invariant of the corresponding auxiliary object data, together with, to every morphism $f: O \rightarrow O'$ of $\cS$ and corresponding set of auxiliary data $\cG(f) \subset \cG(O_1) \times \cG(O_2)$ a map $$F_f(\cD_1, \cD_2): F(\cD_1) \rightarrow F(\cD_2)$$ for every $(\cD_1, \cD_2) \in \cO(f)$. 

We also require that for every pair of edges $(e_1, e_2)$ in $\cE(f)$ which connect $(\cD_1, \cD_2)$ and $(\cD_1', \cD_2')$ in $\cV(f)$, the diagram
$$\begin{tikzcd}F(\cD_1) \arrow[r, "F_f(\cD_1{,}\cD_2)"] \arrow[d, swap, "F(e_1)"]& F(\cD_2) \arrow[d, "F(e_2)"]\\
F(\cD_1') \arrow[r, "F_f(\cD'_1{,}\cD'_2)"] & F(\cD_2')
\end{tikzcd}$$
commutes in $\cC$.
\end{defn}

We may also take preferred triangles into account:
\begin{defn}Let $\cS$ be a category with preferred triangles, and $\cG = \cG(\cS)$ be a set of auxiliary data for $\cS$ with triangles. For any triangulated category $\cC$, a \emph{weak triangulated $\cG$ invariant of $\cS$ (valued in $\cC$)} is a weak invariant of the category $\cS$ valued in $\cC$, together with extra conditions for any preferred triangle $\Delta = (O_i, f_i)_{i=0}^{i=2}$ in $\cS$ and corresponding set of auxiliary data $\cG(\Delta) \subset \cG(O_0) \times \cG(O_1) \times \cG(O_2)$.

Namely, for each $(\cD_0, \cD_1, \cD_2) \in \cG(\Delta)$ we require that the triangle $$(F(\cD_i), F_{f_i}(\cD_i, \cD_{i+1}))_{i=0}^{i=3}$$ is distinguished in $\cC$.
\end{defn}
% Note that as the maps $F(p_i)$ are isomorphisms, these latter morphisms $(F(p_i))^{i=3}_{i=0}$ are, in fact, isomorphisms.

\subsubsection{Strong invariants and transitive systems}
Working with weak invariants of objects and isomorphism classes precludes any useful or interesting morphisms in an invariant manner --- there are no interesting morphisms between isomorphsim classes. To take this extra structure into account, we have to be more stringent:
\begin{defn}A \emph{strong $\cG$ invariant of $\ob(\cS)$ (valued in $\cC$)} is a weak $\cG$ invariant of $\ob(\cS)$ (valued in $\cC$) for which, for any object $O \in \ob(\cS)$, any data $\cD_1, \cD_2 \in \cV(O)$ and any two choices of path $p, p': \cD_1 \rightarrow \cD_2$, the corresponding morphisms $F(p)$ and $F(p')$ are the same.

More specifically, if $\cS$ is a category (with preferred triangles), a \emph{strong (triangulated) $\cG$ invariant of $\cS$ (valued in $\cC$)} is a weak (triangulated) $\cG$ invariant of $\cS$ (valued in $\cC$) for which the corresponding weak $\cG$ invariant of $\ob(\cS)$ is a strong $\cG$ invariant of $\ob(\cS)$.
\end{defn}

A strong invariant of a category $\cS$ valued in $\cC$ is supposed to give a functor $\cS \rightarrow \cC$. The key algebraic set-up  (cf. \cite{jtz:naturality}, \cite{bs:naturalitymonopole}) is that of a \emph{transitive system}.
\begin{defn}[{\cite[Definition 1.1]{jtz:naturality}}]For any category $\cC$, a \emph{transitive system in $\cC$} is
\begin{enumerate}\item a set $M$, and for every $\alpha \in M$ an object $O_\alpha$
\item for every pair $(\alpha, \beta) \in M \times M$ an isomorphism $\pi^\alpha_\beta: O_\alpha \rightarrow O_\beta$ such that $\pi^\beta_\gamma \circ \pi^\alpha_\beta = \pi_\gamma^\alpha$ for every $\alpha, \beta, \gamma \in M$.
\end{enumerate}
\end{defn}

\begin{rmk}In the case where $\cC$ is a \emph{concrete} category, i.e. its objects are sets and morphisms are maps of sets, one can define an honest element of $\cC$ associated to a transitive system, by letting $O(M, \{\pi_\beta^\alpha\})$ be the set of elements $g \in \Pi_{\alpha \in M} O_\alpha$ such that $\pi^\alpha_\beta(g(\alpha)) = g(\beta)$ for every $\alpha, \beta \in M$.
\end{rmk}

\begin{defn}[{\cite[Definition 1.1]{bs:naturalitymonopole}}]A \emph{morphism} of transitive systems in $\cC$ from $(M, \{O_\alpha\}, \{\pi_\beta^\alpha\})$ to $(M', \{O'_\gamma\}, \{\eta_\delta^\gamma\})$ is, for every $(\alpha, \gamma) \in M \times M'$ a morphism $f^\alpha_\gamma: O_\alpha \rightarrow O'_\gamma$ such that $f^\beta_\delta \circ \pi^\alpha_\beta = \pi_\delta^\gamma f^\alpha_\gamma$ for every $\alpha, \beta \in M$ and $\gamma, \delta \in M'$.
\end{defn}

It's worth noting here that, as in \cite[Remark 2.6]{bs:contactsuturedmonopole}, any collection of equivalence classes of morphisms $\{f^\alpha_\gamma\}$ where the indices belong to a nonempty subset of indexing sets $M \times M'$ can be uniquely completed to a morphism of transitive systems so long as they satisfy the properties a morphism has to where they are defined.

\begin{rmk}
When $\cC$ is a concrete category, a morphism of transitive systems $$\{f_\gamma^\alpha\}: (M, \{O_\alpha\}, \{\pi_\beta^\alpha\}) \rightarrow (M', \{O'_\gamma\}, \{\eta_\delta^\gamma\})$$ defines an honest morphism of corresponding objects $$f: O(M, \{O_\alpha\}, \{\pi_\beta^\alpha\}) \rightarrow O(M', \{O'_\gamma\}, \{\eta_\delta^\gamma\}).$$
\end{rmk}

These two notions give us a category, $\csys$; when $\cC$ is a concrete category this comes with a functor to the category $\cC$. When $\cC$ is triangulated, $\csys$ inherits a triangulated structure by declaring the distinguished triangles in $\csys$ to be those where some level-wise triangle is distinguished in $\cC$.

\begin{defn}Let $F$ be a strong (triangulated) $\cG$ invariant of $\cS$ valued in $\cC$. We define a corresponding functor $\cF: \cS \rightarrow \csys$, by dictating
\begin{enumerate}
\item for any $O$ in $\ob(\cS)$ we let $\cF(O)$ be the transitive system with objects given by the $F(\cD)$ for $\cD \in \cG(O)$ and the maps $\pi^{F(\cD_1)}_{F(\cD_2)}$ given by $F(p)$ for any path $p$ connecting $\cD_1$ and $\cD_2$; and
\item for any morphism $f: O_1 \rightarrow O_2$ of $\cS$ we let $\cF(f)$ be the morphism of transitive systems obtained by extending the maps $F(\cD_1, \cD_2)$ for any $(\cD_1 , \cD_2) \in \cG(f)$ to the whole of the transitive systems $\cF(O_1), \cF(O_2)$.\end{enumerate}
\end{defn}

A straightforward consequence of the set-up is:
\begin{prop}Let $\cS$ have preferred triangles, and $\cC$ be triangulated. Any strong triangulated $\cG$ invariant of $\cS$ valued in $\cC$ sends preferred triangles to distinguished triangles.
\end{prop}

In the case where $\cC$ is a concrete category, by post-composing with the functor from $\csys$ to $\cC$ we obtain:
\begin{prop}Let $\cS$ be a category. Any strong $\cG$ invariant of $\cS$ valued in a concrete category $\cC$ gives a functor from $\cS$ to $\cC$.
\end{prop}

\subsection{Homotopy invariants}
Usually, when working with the above set-up, we deal with weak/strong $\cG$ invariants which are valued in the derived $\DC$ or homotopy category $\KC$ associated to a dg category $\cC$. The construction of $F$ in this context usually factors through constructions in $\cC$ ---  then taking their images under the corresponding functor to the derived or homotopy category.

It's usually \emph{necessary} to pass to the derived or homotopy category because we don't have strong enough naturality properties of morphisms: for instance, where the maps $F(p_i)$ should commute with maps $F(\cD_1, \cD_2)$ (cf. Definition \ref{weakinvariant}) they only {homotopy commute}, or where we desire that maps $F(e)$ are isomorphisms they are only quasi-isomorphisms or homotopy equivalences.

It's often \emph{desirable} to pass to the derived or homotopy category because the category $\cC$ is not in and of itself triangulated --- whereas the latter are.

Nevertheless, passing to the derived category or homotopy category comes with its own drawback: morphisms are no longer concrete maps of sets, but equivalence classes of such maps. We also have to deal with one of the painful details of triangulated categories in general: that the mapping cone of a morphism is only well-defined up to isomorphism --- this is particularly apparent in the case of the derived or homotopy category. In what follows we will need to consider morphisms into and out of mapping cones: this only makes sense in the concrete setting where the cone is an actual object of $\cC$, rather than an isomorphism class. This essentially forces us to work in the honest dg category where we have such a notion.

With this in mind, we `lift' the above definitions to constructions in a dg category in such a way that passing to the derived category or homotopy category gives a weak/strong invariant as above.

\begin{defn}Let $\cS$ be a category, and $\cG$ be a set of auxiliary data for $\cS$ (with triangles). For a dg category $\cC$, a \emph{weak (resp. strong) $\cG$ homotopy invariant of $\ob(\cS)$ (valued in $\cC$)} is a functor of graphs $F: \cG(\cS) \rightarrow \cC$, where if we denote by $KF: \cG(\cS) \rightarrow \KC$ the data defined by taking $KF(\cD) := F(\cD)$ for all $\cD \in \cV(\cS)$ and $KF(e) = [F(e)]$ for all $e \in \cE(\cS)$, we require that $KF$ is a weak (resp. strong) $\cG$ invariant of $\ob(\cS)$ valued in $\KC$.

Similarly, a \emph{weak (resp. strong) $\cG$ homotopy invariant of $\cS$ (valued in $\cC$)} is a weak homotopy invariant of $\ob(\cS)$ valued in $\cC$ together with, for every morphism $f$ of $\cS$ a subgraph $\cG(f)$ of $\cG(O) \times \cG(O')$ --- with vertices $\cV(f) \subset \cV(O) \times \cV(O')$ and edges $\cE(f) \subset \cE(O) \times \cE(O')$ such that, denoting by $KF(\cG(f))$ the image $KF(\cG(f)) = [F(\cG(f))]$, we require that this additional data constitutes a weak (resp. strong)  $\cG$ invariant of $\cS$ valued in $\KC$.

Finally, a \emph{weak (resp. strong) $\cG$ triangulated homotopy invariant of $\cS$ (valued in $\cC$)} is a weak homotopy invariant of $\cS$ valued in $\cC$ such that for every preferred triangle $\Delta = (O_i, f_i)^{i=2}_{i=0}$ in $\cS$ and $(\cD_0, \cD_1, \cD_2) \in \cO(\Delta)$, the image of the corresponding triangle 
$$(F(\cD_i), F_{f_i}(\cD_i, \cD_{i+1}))_{i=0}^{i=3}$$ 
in $\KC$ is distinguished.
\end{defn}

At the expense of brevity, it's worth belabouring what this means. For a weak $\cG$ triangulated homotopy invariant valued in $\cC$, we require that for any edge $e \in \cE(\cS)$ the morphisms $F(e)$ are homotopy equivalences in $\cC$. Any morphism $F_f(\cD_1, \cD_2)$ is required to homotopy commute with $F(e_1), F(e_2)$ for any pair $(e_1, e_2) \in \cE(f)$. Finally, for any $(\cD_0, \cD_1, \cD_2) \in \cO(\Delta)$ there needs to exist a chain homotopy equivalence from $\cone(F_{f_0}(\cD_0, \cD_1))$ to $F(\cD_2)$.

\subsubsection{Witnesses}

If one wants to show that a given $F$ is a weak homotopy invariant, one has to demonstrate the existence of the homotopies and homotopy equivalences required above. Oftentimes one actually constructs \emph{explicit} maps depending on the data $\cG(\cS)$, that witness this. Given the way these are often constructed, it is a reasonable question to ask how natural these witnesses are in and of themselves.
\begin{defn}Let $F$ be a weak $\cG$ homotopy invariant of $\cS$ (valued in $\cC$). We will say that a set of \emph{witnesses} for $F$ is the further data of, for every triple $(O_1, O_2, f)$ where $O_1, O_2 \in \ob(\cS)$ and $f: O_1 \rightarrow O_2$, any pairs $(\cD_1, \cD_2)$ and $(\cD_1', \cD_2')$ in $\cV(f)$ and any pair of paths $(p_1, p_2) \in \cE(f)$ connecting $(\cD_1, \cD_2)$ to $(\cD_1', \cD_2')$, a choice of explicit homotopy $F_f(p_1, p_2)$ from $F(\cD_1)$ to $F(\cD_2')$ such that:
\begin{enumerate}
\item   the diagram 
$$\begin{tikzcd}[column sep=60pt]F(\cD_1) \arrow[r, "F_{f}(\cD_1{,}\cD_2)"] \arrow[d, swap, "F(p_1)"]& F(\cD_2) \arrow[d, "F(p_2)"]\\
F(\cD_1') \arrow[r, swap, "F_f(\cD'_1{,}\cD'_2)"] & F(\cD_2')
\end{tikzcd}$$
homotopy commutes via the homotopy $F_f(p_1, p_2)$, and

\item the induced morphism $\bF_f(p_1, p_2)$ defined by $$(F(p_1), F(p_2) + F(p_1, p_2)): \cone(F_f(\cD_1, \cD_2)) \rightarrow \cone(F_f(\cD_1', \cD_2'))$$
is a chain homotopy equivalence.
\end{enumerate}
We sometimes say that a pair of weak $\cG$ homotopy invariant of $\cS$ valued in $\cC$ together with a set of witnesses is a \emph{witnessed weak $\cG$ homotopy invariant of $\cS$, valued in $\cC$}.

A \emph{witnessed strong $\cG$ homotopy invariant of $\cS$, valued in $\cC$} is a witnessed weak $\cG$ homotopy invariant of $\cS$ valued in $\cC$ such that the homotopy class of the morphisms $\bF_f(p_1, p_2)$ are independent of choice of $(p_1, p_2) \in \cE(f)$.\end{defn}

One way to view this data is as an invariant of the \emph{arrow category} $\ob(\cS(2))$. This is the category whose objects are triples $(O_1, O_2; f)$ of objects $O_1, O_2$ in $\cS$ and a morphism $f: O_1 \rightarrow O_2$, and where the morphisms between triples $(O_1, O_2; f)$ and $(O_1', O_2'; f')$ are given by pairs $(\phi_1, \phi_2)$ where $\phi_i: O_i \rightarrow O'_i$ are morphisms satisfying $f' \circ \phi_1 = \phi_2 \circ f$.

\begin{prop}Let $F$ be a weak (resp. strong) $\cG$ homotopy invariant of $\cS$ valued in a dg category $\cC$. If $F$ is witnessed then $F$ induces a weak (resp. strong) homotopy invariant of $\ob(\cS(2))$ valued in $\cC$, where $\cS(2)$ is the arrow category of $\cS$.
\end{prop}
\begin{proof}For any triple $(O_1, O_2, f)$ in $\ob(\cS(2))$ we associate the data $\cG(O_1, O_2,f):= \cG(f)$, and define $F(O_1, O_2, f)$ to be $\cone(F(f))$, and $F(p_1, p_2) = \bF_f(p_1, p_2)$.
\end{proof}

Of course, we can address triangles as well:
\begin{defn}Let $F$ be a weak $\cG$ triangulated homotopy invariant of $\cS$ (valued in $\cC$). A set of \emph{witnesses} for $F$ is a set of witnesses for the underlying weak homotopy invariant $F$ together with, for every preferred triangle $\Delta = (O_i, f_i)_{i=0}^{i=2}$ in $\cS$ and any triple $(\cD_0, \cD_1, \cD_2)$ in $\cG(\Delta)$, a choice of homotopy equivalence $$F_\Delta(\cD_0, \cD_1, \cD_2): \cone(F_{f_0}(\cD_0, \cD_1)) \rightarrow F(\cD_2).$$

We also ask for, for any two triples $(\cD_0, \cD_1, \cD_2)$ and $(\cD_0', \cD_1', \cD_2')$ in $\cV(\Delta)$ and paths $(p_0, p_1, p_2)$ in $\cE(\Delta)$, a homotopy $F_\Delta(p_0, p_1, p_2)$ from $\cone(F_{f_0}(\cD_0, \cD_1))$ to $F(\cD_2')$ such that 
\begin{enumerate}
\item the diagram
$$\begin{tikzcd}[column sep=60pt]\cone(F_{f_0}(\cD_0, \cD_1))  \arrow[r, "F_\Delta(\cD_0{,} \cD_1{,} \cD_2)"] \arrow[d, swap, "\bF_{f_0}(p_0{,} p_1)"]& F(\cD_2) \arrow[d, "F(p_2)"]\\
\cone(F_{f_0}(\cD'_0, \cD'_1)) \arrow[r, "F_\Delta(\cD'_0{,} \cD'_1{,} \cD'_2)"] & F(\cD_2')
\end{tikzcd}$$
homotopy commutes via the homotopy $F_\Delta(p_0, p_1, p_2)$, and
\item the induced map $$\bF_\Delta(p_0, p_1, p_2): \cone(F_\Delta(\cD_0{,} \cD_1{,} \cD_2)) \rightarrow  \cone(F_\Delta(\cD'_0{,} \cD'_1{,} \cD'_2))$$ defined by $$(\bF(p_0{,} p_1), F(p_2) + F_\Delta(p_1, p_2, p_3))$$
is a chain homotopy equivalence.
\end{enumerate}

A set of \emph{biwitnesses} for $F$ is a set of witnesses for $F$, together with a choice of homotopy inverse for $F_\Delta(\cD_0, \cD_1, \cD_2)$, written as
$$F_\nabla(\cD_0, \cD_1, \cD_2): F(\cD_2) \rightarrow \cone(F_{f_0}(\cD_0, \cD_1)),$$ along with, for
any two triples $(\cD_0, \cD_1, \cD_2)$ and $(\cD_0', \cD_1', \cD_2')$ in $\cV(\Delta)$ and paths $(p_0, p_1, p_2)$ in $\cE(\Delta)$, a homotopy $F_\nabla(p_0, p_1, p_2)$ from $F(\cD_1)$ to  $\cone(F_{f_0}(\cD_0', \cD_1'))$ such that 
\begin{enumerate}
\item the diagram
$$\begin{tikzcd}[column sep=60pt] F(\cD_2)\arrow[r, "F_\nabla(\cD_0{,} \cD_1{,} \cD_2)"] \arrow[d, swap, "F(p_2)"]&  \cone(F_{f_0}(\cD_0, \cD_1))  \arrow[d, " \bF_{f_0}(p_0{,} p_1)"]\\
F(\cD_2') \arrow[r, "F_\nabla(\cD'_0{,} \cD'_1{,} \cD'_2)"] & \cone(F_{f_0}(\cD'_0, \cD'_1)) 
\end{tikzcd}$$
homotopy commutes via the homotopy $F_\nabla(p_0, p_1, p_2)$, and
\item the induced map $$\bF_\nabla(p_0, p_1, p_2): \cone(F_\nabla(\cD_0{,} \cD_1{,} \cD_2)) \rightarrow  \cone(F_\nabla(\cD'_0{,} \cD'_1{,} \cD'_2))$$ defined by $$(F(p_2),   \bF(p_0{,} p_1)+ F_\nabla(p_1, p_2, p_3))$$
is a chain homotopy equivalence.
\end{enumerate}

We sometimes say that a pair a of weak $\cG$ homotopy invariant of $\cS$ valued in $\cC$ and a set of (bi)witnesses is a \emph{(bi)witnessed weak $\cG$ homotopy invariant of $\cS$, valued in $\cC$}.

A \emph{witnessed (resp. biwitnessed) strong $\cG$ homotopy invariant of $\cS$, valued in $\cC$} is a (bi)witnessed weak $\cG$ homotopy invariant of $\cS$ valued in $\cC$ such that the homotopy classes of the morphisms $\bF_\Delta(p_0, p_1, p_2)$ are independent of choice of $(p_0, p_1, p_2) \in \cE(f)$. (Resp. this is also true for the homotopy classes of the morphisms $\bF_\nabla(p_0, p_1, p_2)$.)
\end{defn}
\subsubsection{Bordered-sutured categories and carriage returns}\label{carriagereturns}
If we are concerned with the category $\cS = \contdiffborsut$, there is a weakening of the set of witnesses for a weak triangulated invariant which we will care about.

As motivation and, as it happens, what appears to be the most common situation in the `wild', it is illuminating to first consider properties of explicit weak triangulated $\cG$ homotopy invariants of $\contdiffborsut$.

A particularly interesting property of $\contdiffborsut$ is that some triangles come in conjugate pairs. A (strongly) witnessed weak triangulated $\cG$ homotopy invariant of $\contdiffborsut$ valued in a dg category $\cC$ associates to these two \emph{a priori} different sets of auxiliary data. We could ask for some compatibility for these:
\begin{defn}Let $F$ be a biwitnessed weak triangulated $\cG$ homotopy invariant of $\contdiffborsut$ (valued in $\cC$). We will say that $F$ \emph{admits carriages} if, for every pair of conjugate triangles  $\Delta_f = (O_i, f_i)_{i=0}^{i=3}$ and  $\Delta_h = (O_i, h_i)_{i=0}^{i=3}$  of $\contdiffborsut$ there are sets $(\cD_0, \cD_1, \cD_2) \in \cV_{\Delta_f}(O_0, O_1, O_2))$ and $(\cD'_0, \cD'_1, \cD'_2) \in \cV_{\Delta_h}(O_0, O_1, O_2))$ such that $F(\cD_i) = F(\cD'_i)$.
\end{defn}

In this case, there are chain homotopy equivalences $$F(\cD_0, \cD_1; \cD_0', \cD'_1): \cone(F_{f}(\cD_0, \cD_1)) \rightarrow \cone(F_{h}(\cD'_0, \cD'_1))$$ defined by the formula
$$F(\cD_0, \cD_1; \cD_0', \cD'_1): = F_{\nabla_h}(\cD'_0, \cD'_1, \cD'_2) \circ F_{\Delta_f}(\cD_0, \cD_1, \cD_2).$$

These satisfy an invariance property, which we encode in a definition:
\begin{defn}A \emph{weak $\cG$ homotopy invariant of $\contdiffborsut$ with carriage returns (valued in $\cC$)} is a weak $\cG$ homotopy invariant of $\contdiffborsut$ (valued in $\cC$) together with, for every pair $(Y, f) \in \ob(\diffborsuttor_2^\perp)$, a choice of nonempty connected subgraph $$\cG(Y, f(Y)) \subset \{((\cD_1, \cD_2),(\cD'_1, \cD'_2)) \in \cG(f) \times \cG(h): F(\cD_i) = F(\cD'_i) \mbox{ for } i = 0, 1\}$$
and an assignment of, for any $((\cD_0, \cD_1), (\cD_0', \cD_1')) \in \cV(Y, f(Y))$, a morphism
$$F(\cD_0, \cD_1; \cD_0', \cD'_1): \cone(F_{f}(\cD_0, \cD_1)) \rightarrow \cone(F_{h}(\cD'_0, \cD'_1)),$$
and, for any quadruple of paths $((p_0, p_1); (p_0', p_1'))$ of $\cE(Y, f(Y))$ connecting $((\cD_0, \cD_1);(\cD'_0, \cD'_1))$ to $((\cD_0'', \cD_1'');(\cD_0''', \cD_1'''))$, a homotopy $$F(p_0, p_1; p_0', p_1'):\cone(F_f(\cD_0, \cD_1))\rightarrow \cone(F_h(\cD_0''', \cD_1''')),$$ such that

\begin{enumerate}\item the diagram
$$\begin{tikzcd}[column sep=80pt]\cone(F_f(\cD_0, \cD_1))  \arrow[r, "F(\cD_0{,} \cD_1; \cD_0'{,} \cD'_1)"] \arrow[d, swap, "\bF_{\Delta_f}(p_0{,} p_1)"]& \cone(F_h(\cD'_0, \cD'_1))\arrow[d, "\bF_{\Delta_h}(p_1' {,} p_2')"]\\
\cone(F_f(\cD''_0, \cD''_1))  \arrow[r, "F(\cD''_0{,} \cD''_1; \cD'''_0{,} \cD'''_1)"] & \cone(F_h(\cD'''_0, \cD'''_1))
\end{tikzcd}$$
homotopy commutes via the homotopy $F(p_0, p_1; p_0', p_1')$, and 
\item the induced map
$$\bF(p_0, p_1; p'_0, p_1'): \cone(F(\cD_0{,} \cD_1; \cD_0'{,} \cD'_1)) \rightarrow  \cone(F(\cD''_0{,} \cD''_1; \cD'''_0{,} \cD'''_1))$$
defined by $$(\bF_{\Delta_f}(p_0{,} p_1), \bF_{\Delta_h}(p'_0{,} p'_1), + F(p_0, p_1; p_0', p_1'))$$
is a chain homotopy equivalence.
\end{enumerate}
\end{defn}

\begin{rmk}\label{carriagedeparturesfromtriangulated}The situation in a biwitnessed weakly triangulated invariant which admits carriages is symmetric, and we also obtain 
$$F(\cD_0, \cD_1; \cD_0', \cD'_1): \cone(F_{h}(\cD_0, \cD_1)) \rightarrow \cone(F_{f}(\cD'_0, \cD'_1))$$ defined by the formula
$$F(\cD_0, \cD_1; \cD_0', \cD'_1): = F_{\nabla_f}(\cD'_0, \cD'_1, \cD'_2) \circ F_{\Delta_h}(\cD_0, \cD_1, \cD_2),$$
which is a homotopy inverse to $F(\cD_0, \cD_1; \cD_0', \cD_1')$

One can codify this and the corresponding invariance properties in a similar fashion to define a \emph{weak $\cG$ homotopy invariant of $\contdiffborsut$ with carriages (valued in $\cC$)} --- we avoid belabouring the point here and return to it later with some more notation in hand.
\end{rmk}

\subsection{The category of type-writer structures}
Just as a consistent weak homotopy invariant of $\cS$ gives a weak homotopy invariant of the objects of the arrow category $\cS(2)$, it turns out that the data of a weak homotopy invariant of $\contdiffborsut$  with carriage returns valued in $\cC$ gives a weak homotopy invariant of the category $\diffborsuttor_2^\perp$, valued in a different category.

\begin{defn}\label{typewriterdefn}Let $\cC$ be a dg category. The \emph{category of typewriters associated to $\cC$}, written $\crC$ is the dg category whose objects are tuples $$\bM := (M_0, M_1; D_f, D_h; D_{CR})$$ called \emph{typewriters}, where:
\begin{enumerate}\item $M_0$ and $M_1$ are objects of $\cC$.
\item $D_f, D_h: M_0 \rightarrow M_1$ are morphisms.
\item $D_{CR}: \cone(D_f) \rightarrow \cone(D_h)$ is a morphism, called a \emph{carriage return}.
\end{enumerate}

A \emph{morphism} of typewriters from $\bM = (M_0, M_1; D_f, D_h; D_{CR})$ to $\bM' = (M'_0, M'_1; D'_f, D'_h; D'_{CR})$ is a tuple $\bT = (T_0, T_1; T_f, T_h; T_{CR})$ where:
\begin{enumerate}\item $T_i: M_i \rightarrow M_i'$ are morphisms for $i = 0, 1$;
\item $T_*$ for $* = f, h$ is a homotopy $M_0 \rightarrow M_1'$ such that the diagram
$$\begin{tikzcd}M_0 \arrow[r, "D_*"] \arrow[d, swap, "T_0"]& M_1 \arrow[d, "T_1"]\\
M'_0 \arrow[r, "D'_*"] & M'_1
\end{tikzcd}$$
homotopy commutes via the homotopy $T_*$, inducing maps
$$\bT_*: \cone(D_*) \rightarrow \cone(D'_*)$$ via the formula
$$\bT_*(x,y) :=  (T^0(x), T^1(y) + T^f(x));$$
 and
\item the map $T_{CR}$ is a homotopy $\cone(D_f) \rightarrow \cone(D_h')$ such that the diagram
$$\begin{tikzcd}\cone(D_f) \arrow[r, "D_{CR}"] \arrow[d, swap, "\bT_f"]& \cone(D_h) \arrow[d, "\bT_h"]\\
\cone(D'_f) \arrow[r, "D'_{CR}"] & \cone(D'_h)
\end{tikzcd}$$
homotopy commutes via the homotopy $T_{CR}$, inducing a map
$$\bT_{CR}: \cone(D_{CR}) \rightarrow \cone(D'_{CR})$$
by the formula
$$\bT_{CR}(\bx, \by) = (\bT_f(\bx), \bT_h(\by) + T_{CR}(\bx)).$$
\end{enumerate}

We define the composition of two morphisms $\bT: \bM \rightarrow \bM'$ and $\bT':\bM' \rightarrow \bM''$ to be the unique morphism $\bT' \circ \bT= (U_0, U_1; U_f, U_h; U_{CR})$ such that $U_i = T'_i \circ T_i$, $\bU_* = \bT_*' \circ \bT_*$ for $* = f, h$ and $\bU_{CR} = \bT_{CR}' \circ \bT$.
\end{defn}

\begin{exmp}For any typewriter $\cR = (M_0, M_1; D_f, D_h; D_{CR})$ there is a corresponding \emph{identity morphism} $\II: \cR \rightarrow \cR$, given by $\II = (\id_{M_0}, \id_{M_1}; 0; 0)$.
\end{exmp}

There is an appropriate notion of a \emph{homotopy} of morphisms of carriage returns:
\begin{defn}Let $\bM$ and $\bM'$ be typewriters as above, and let $\bT = (T_0, T_1; T_f, T_h; T_{CR})$ and  $\bT' = (T_0', T_1'; T_f', T_h'; T_{CR}')$ be morphisms $\bT, \bT': \bM \rightarrow \bM'$. A \emph{homotopy from $\bT$ to $\bT'$} is a tuple $\bH = (H_0, H_1; H_f, H_h; H_{CR})$ where:
\begin{enumerate}\item $H_i: M_i \rightarrow M_i'$ is a homotopy between $T_i$ and $T'_i$ for $i = 0, 1$;
\item $H_*$ for $* = f, h$ is a map $M_0 \rightarrow M_1'$ which induces a homotopy
$$\bH_*: \cone(D_*) \rightarrow \cone(D'_*)$$ from $\bT_*$ to $\bT_*'$ via the formula
$$\bH_*(x,y) :=  (H_0(x), H_1(y) + H_*(x));$$
 and
\item $H_{CR}$ is a map $\cone(D_f) \rightarrow \cone(D_h')$ such that the induced map $$\bH_{CR}: \cone(T_{CR}) \rightarrow \cone(T_{CR}')$$
given by the formula
$$\bH_{CR}(\bx, \by) = (\bH_f(\bx), \bH_h(\by) + H_{CR}(\bx))$$
is a homotopy.
\end{enumerate}
\end{defn}

This begets the definition:
\begin{defn}A morphism of typewriters $\bT: \bM \rightarrow \bM'$ is a \emph{homotopy equivalence} if it has a \emph{homotopy inverse} $\bT'$ such that $\bT' \circ \bT$ is homotopic to the identity map of $\bM$.

Two typewriters are \emph{homotopy equivalent} if there exists a homotopy equivalence between them.
\end{defn}

With this in mind:
\begin{prop}\label{carriagereturnsistypewriters}A weak (resp. strong) $\cG$ homotopy invariant for $\contdiffborsut$ with carriage returns valued in $\cC$ yields a weak (resp. strong) homotopy invariant for $\diffborsuttor_2^\perp$, valued in $\crC$.
\end{prop}
\begin{proof}This is a straightforward, if somewhat notationally troubling unwrapping of the definitions.
\end{proof}

We can also work with biwitnessed weak triangulated invariants:
\begin{prop}\label{triangulatedistypewriters}Any biwitnessed weak (resp. strong) triangulated $\cG$ homotopy invariant of $\contdiffborsut$ valued in $\cC$ which admits returns yields a weak (resp. strong) homotopy invariant of $\diffborsuttor_2^\perp$, valued in $\crC$.
\end{prop}
\begin{proof}This follows from the discussion in Remark \ref{carriagedeparturesfromtriangulated} together with Proposition \ref{carriagereturnsistypewriters}.
\end{proof}

\subsubsection{Partially extendable typewriters and carriage departures}
In Remark \ref{carriagedeparturesfromtriangulated}  we observed that a biwitnessed weak triangulated invariant of $\contdiffborsut$ carries more information than a weak homotopy invariant with carriage returns. Here, we sketch the construction of the data needed to capture this behaviour.

The appropriate notion is:
\begin{defn}Let $\cC$ be a dg category. The \emph{category of partially extendable typewriters associated to $\cC$}, written $\crpC$, is the subcategory of typewriters whose objects $$\bM= (M_0, M_1; D_f, D_h; D_{CR})$$ satisfy that $D_{CR}$ is a homotopy equivalence.

We often call a homotopy inverse to $D_{CR}$ a \emph{carriage departure}.
\end{defn}

This allows for a refinement of Proposition \ref{triangulatedistypewriters}:
\begin{prop}\label{triangulatedistypewriters}Any biwitnessed weak triangulated $\cG$ invariant of $\contdiffborsut$ valued in $\cC$ which admits returns yields a weak homotopy invariant of $\diffborsuttor_2^\perp$, valued in $\crpC$.
\end{prop}
\begin{proof}This follows closely along the lines of the proof of Proposition \ref{carriagereturnsistypewriters}.
\end{proof}

\begin{rmk}One could strengthen this further and define a \emph{category of partially extended typewriters associated to $\cC$}, where the objects are partially extendable typewriters together with a \emph{specific choice} of carriage departure $D_{CD}$ --- and where morphisms need to respect this in a similar manner to the way they respect the carriage return $D_{CR}$.

It follows by Remark \ref{carriagedeparturesfromtriangulated}  that a biwitnessed weak triangulated $\cG$ invariant of $\contdiffborsut$ valued in $\cC$ which admits returns gives a weak homotopy invariant of $\diffborsuttor^\perp_2$, valued in this category of partially extended typewriters --- i.e. it remembers specific homotopy inverses for the carriage return.
\end{rmk}

\section{Typewriter structures and type $D$ structures}
The aim of this section is rather simple: we show that a category of typewriter structures and a category of type $D$ structures are the same. We restrict our attention to the category of typewriters in a category of type $D$ structures, $\amod$.

\begin{prop}\label{typewritersaretyped}Let $\bM = (M_0, M_1; D_f, D_g, D_{CR})$ be a typewriter in $\amod$. Then there is a corresponding $DD$ bimodule $\dify(\bM) \in \ob(\amodt)$.
\end{prop}
\begin{proof}\label{functoronobjects}
Splitting the morphism $D_{CR}: \cone(D_f) \rightarrow \cone(D_g)$ into components, we see that this gives four maps $D_{g}, D_{fg}, D_{gh}, D_{fgh}$ such that for a pair $(x, y) \in \cone(D_f)$ we have that $$D_{CR}(x,y) := (D_{g}(y) + D_{fg}(x), D_{gh}(y) + D_{fgh}(x)).$$

We claim that the following expressions all evaluate to zero:
\begin{alignat*}{3}
D_\emptyset \circ D_f  + D_f \circ D_\emptyset &  \qquad D_\emptyset \circ D_g  + D_g \circ D_\emptyset  &\qquad  D_\emptyset \circ D_h  + D_h \circ D_\emptyset \end{alignat*}
\begin{alignat*}{2}
  D_{fg}\circ D_\emptyset + D_\emptyset \circ D_{fg} + D_g \circ D_f  &\qquad  D_{gh} \circ D_\emptyset + D_\emptyset \circ D_{gh} + D_h \circ D_g  
\end{alignat*}
\begin{alignat*}{1} D_{fgh}\circ D_\emptyset + D_\emptyset \circ D_{fgh} + D_{gh} \circ D_f +D_h \circ D_{fg} \end{alignat*}

Indeed, the expressions in the first line say that $D_f, D_g, D_h$ are all morphisms in $\cD$. The first and the last of these come for free by definition; we defer the middle for a moment.

As $D_{CR}(x,y)$ is a chain map, we must have that $D_{CR} \circ \partial_{\cone(D_f)} + \partial_{\cone(D_h)} \circ D_{CR} = 0$. Explicitly resolving this into components gives
$$D_{CR} \circ \partial_{\cone(D_f)}(x, y) = (D_g \circ D_f(x)) + D_g\circ D_\emptyset'(y) + D_{fg} \circ D_\emptyset(x), D_{gh} \circ D_f(x) + D_{gh} \circ D_\emptyset'(y) + D_{fgh} \circ D_\emptyset(x))$$
and
$$\partial_{\cone(D_h)} \circ D_{CR}(x,y) = (D_\emptyset \circ D_{g}(y) + D_\emptyset \circ D_{fg}(x), D_h \circ D_{g}(y) + D_h \circ D_{fg}(x) +D_\emptyset' \circ D_{gh}(y) + D_\emptyset' \circ D_{fgh}(x)).$$
Evaluating at $(x, 0)$ and comparing terms gives the first of the middle line of expressions and the final expression. Similarly, evaluating at $(0, y)$ gives the middle of the first line of expressions and the second of the middle line of expressions.

Therefore, by Proposition \ref{coefftodd}, if we let $\dify(\bM)$ be given by the direct sum $M_0 \oplus M_1$ and define a map $\delta^1 : \dify(\bM) \rightarrow \cA \otimes \dify(\bM) \otimes \cT$ by the formula
$$\delta^1 := (D_\emptyset, D_\emptyset') \otimes \bone + \sum_{t \in \cT'} D_{t} \otimes t,$$
the pair $(\dify(\bM), \delta^1)$ is a type $DD$ bimodule over $(\cA, \cT)$.
\end{proof}

Similarly, we can define $\dify$ on morphisms:
\begin{prop}
Let $\bT: \bM \rightarrow \bM'$ be a morphism of typewriters, $\bT:=(T_0, T_1; T_f, T_h; T_{CR})$. Then there is a corresponding morphism of $DD$ bimodules $\dify(\bT): \dify(\bM) \rightarrow \dify(\bM)$. Moreover, given two morphisms of typewriters $\bT: \bM\rightarrow \bM'$ and $\bT': \bM' \rightarrow \bM''$, we have that $\dify(\bT' \circ \bT) = \dify(\bT') \circ \dify(\bT)$.
\end{prop}
\begin{proof}
This follows similarly to Proposition \ref{functoronobjects}; resolving the map $$T_{CR}: \cone(D_f) \rightarrow \cone(D'_h)$$ into components, we write:
$$T_{CR}(x,y) :=  (T_{g}(y) + T_{fg}(x), T_{gh}(y) + T_{fgh}(x)),$$
and define a map  $\psi^1: \dify(\bM) \rightarrow \cA \otimes \dify(\bM') \otimes  \cT' $ by the formula
$$\psi^1 :=  (T_0, T_1) \otimes \bone + \sum_{t \in \cT'} T_{t} \otimes t .$$

One can verify that this is, indeed, a morphism of type $DD$ bimodules. Moreover, a (somewhat arduous) calculation shows the composition property in the statement also holds.
\end{proof}

\begin{exmp}\label{idsenttoid}Let $\bM$ be a typewriter, and $\II: \bM \rightarrow \bM$ be the identity morphism, which we recall is given by $(\id_{M_0}, \id_{M_1}; 0, 0; 0)$. We have that the corresponding morphism $\dify(\II): \dify(\bM) \rightarrow \dify(\bM)$ is given by $\dify(\II)(\bx) = \bx \otimes \bone$, i.e. $\dify$ sends the identity morphism to the identity morphism.
\end{exmp}

Finally, we consider homotopies:
\begin{prop}\label{homotopiestohomotopies}
Let $\bH: \bT \rightarrow \bT'$ be a homotopy of morphisms of typewriters in $\amod$, $\bT, \bT':  \bM \rightarrow \bM'$, given by $\bH = (H_0, H_1; H_f, H_h; H_{CR})$. Then there is a corresponding homotopy of $DD$ bimodules $\dify(\bH) : \dify(\bT) \rightarrow \dify(\bT')$.

Moreover, if $\bT: \bM \rightarrow \bM'$ is a homotopy equivalence then so is $\dify(\bT)$.
\end{prop}
\begin{proof}
Again, this follows very similarly to Proposition \ref{functoronobjects}.  Resolving the map $\bH_{CR}$ into components gives
$$H_{CR}(x,y) :=  (H_{g}(y) + H_{fg}(x), H_{gh}(y) + H_{fgh}(x)).$$
and we define a map  $\bH^1: \dify(\bM) \rightarrow \cA \otimes \dify(\bM')\otimes \cT$ by the formula
$$\bH^1 := (H_0, H_1) \otimes \bone + \sum_{t \in \cT}  H_{t} \otimes t,$$
which one may verify is a homotopy of $DD$ bimodules.

The final sentence of the statement follows from this together with Example \ref{idsenttoid}.
\end{proof}

\begin{thm}\label{functorequivalence}The functor $\dify$ is an equivalence of categories, which descends to the homotopy categories of both.
\end{thm}
\begin{proof}As constructed $\dify$ is clearly full, faithful and essentially surjective. It follows from Proposition \ref{homotopiestohomotopies} that this descends to the homotopy categories.
\end{proof}

This leads us to the following:
\begin{cor}\label{carriagegivestyped}Let $F$ be a weak (resp. strong) witnessed $\cG$ homotopy invariant of $\contdiffborsut$ valued in $\amod$, with carriage returns. Then there is a corresponding weak (resp. strong) homotopy $\cG$ invariant of $\diffborsuttor^\perp_2$ valued in $\amodt$ which we denote by $BFD$.\end{cor}
\begin{proof}This follows readily from Proposition \ref{carriagereturnsistypewriters} together with Theorem \ref{functorequivalence}.
\end{proof}

\subsection{Semi-extendibility}
We can play a more elaborate game with the category of semi-extendable typewriters.

\begin{prop}Let $\bM = (M_0, M_1; D_f, D_h; D_{CR})$ be a semi-extendable typewriter valued in $\amod$. Then $\dify(\bM)$ is semi-extendable.
\end{prop}
\begin{proof}Let $D_{CD}$ be a carriage departure for $D_{CR}$. Splitting the morphism $$D_{CD}: \cone(D_h) \rightarrow \cone(D_f)$$ into components, we see that this gives four maps  $D_e, D_{he}, D_{hef}, D_{ef}$, and that these satisfy half of the analogous relations for an semi-extendable type $DD$ bimodule, Definition \ref{semiextendable}.

Using the information that $D_{CD} \circ D_{CR}$ is chain homotopic to the identity via some homotopy $H_{C}: \cone(D_f) \rightarrow \cone(D_f)$, and resolving $H_C$ into components $$(D_{fghe}(x) + D_{ghe}(y), D_{fghef}(x)+D_{ghef}(y))$$
yields two of the corresponding structure relations for these generalised coefficient maps. The other two structure relations are obtained in a similar manner, resolving a homotopy between $D_{CR} \circ D_{CD}$ and the identity.
\end{proof}

\begin{rmk}It is worth remarking upon the appearance of the map $D_{fghef}$ here. Although this is unnecessary in verifying that the corresponding type $D$ structure is semi-extendable, this form of structure map \emph{does} appear in the definition of an extendable type $D$ structure --- and the relation one obtains from the above construction agrees with that which it has to satisfy.\end{rmk}

\begin{rmk}Of course, if we have a semi-extended typewriter then we can say more: to $\bM = (M_0, M_1; D_f, D_h; D_{CR}, D_{CD})$ the above proposition associates an explicit set of semi-extended coefficient maps. Presumably this extends to a functor: the details of defining morphisms and homotopies between semi-extended $DD$ bimodules is rather laborious, so we don't pursue this avenue here.
\end{rmk}

\subsection{2-categories}
\label{typewritersaretwocategory}
In a debt to the introduction, we conclude this section by defining the horizontal composition which realises typewriters as cells in a 2-catgory. This isn't particularly useful as part of the general thrust of the paper, and the reader is invited to skip it.

\begin{prop}Let $\bM = (M_0, M_1; D_f, D_h; D_{CR})$ and $\bM' = (M_1, M_2; D'_f, D'_h; D'_{CR})$ be typewriters. Then there is a well-defined typewriter $\bM \star \bM'$.
\end{prop}
\begin{proof}As in the proof of Proposition \ref{typewritersaretyped}, we can split $D_{CR}$ and $D'_{CR}$ into components,
\begin{align*}D_{CR}(x,y) := (D_{g}(y) + D_{fg}(x), D_{gh}(y) + D_{fgh}(x)).\\
D'_{CR}(x,y) := (D'_{g}(y) + D'_{fg}(x), D'_{gh}(y) + D'_{fgh}(x)).\end{align*}
We define $$\bM \star \bM' = (M_0, M_2; (D \star D')_f, (D \star D')_h; (D \star D')_{CR})$$ by
\begin{align*}(D \star D')_f & := D'_f \circ D_f\\
(D \star D')_h & := D'_h \circ D_h\end{align*}
and setting
$$ (D \star D')_{CR}(x, y) := (D^\star_{g}(y) + D^\star_{fg}(x), D^\star_{gh}(y) + D^\star_{fgh}(x)),$$
where
\begin{align*}D^\star_g &:= D_g \circ D_g'\\
D^\star_{fg} & := D_g \circ D_{fg}' \circ D_f\\
D^\star_{gh} &:= D'_h \circ D_{gh} \circ D'_g\\
D^\star_{fgh} & := D_h' \circ D_{gh} \circ D'_{fg} \circ D_f.
\end{align*}
\end{proof}

\section{Bordered invariants compatible with sutured Floer homology}
\label{borderedsuturedsection}

In this section, we work with the weak invariant $\bsd$ and construct appropriate carriage returns --- ultimately, by Corollary \ref{carriagegivestyped} this gives us a weak invariant of $\diffborsuttor_2^\perp$ which we show to be the same as the corresponding $\bsd$ invariant. The first step to this is to address the case of the thickened torus, and define a corresponding carriage return here. We then use this model computation to extend the construction to deal with the whole category $\contdiffborsut$.

\subsection{The thickened torus}
We spend this section on the example where $Y$ is a thickened torus $T^2 \times I$, endowed with a parameterisation of $T^2 \times \{0\}$ by an arc diagram for the torus. We are concerned with the case where $\gamma$ is a meridian $\mu$ and $f(\gamma)$ is a standard, $0$-framed longitude $\lambda$ respectively, and run the above procedure to produce a corresponding typewriter.

\subsubsection{Type $D$ structure representatives}

We fix a pair of bordered-sutured diagrams $\cD(\mu)$ and $\cD(\lambda)$ representing the bordered-sutured three-manifolds $Y$ and $f(Y)$ respectively --- see Figure \ref{fig:torusdiagrams}. With these chosen, we define $F(\cD(\mu))$ and $F(\cD(\lambda))$  to be the type $D$ structures $\bsd(\cD(\mu))$ and $\bsd(\cD(\lambda))$ respectively.
\begin{figure}
  \begin{minipage}{0.3\linewidth}\centering
    %\scalebox{1}{\includesvg[pretex=\relscale{1}]{heegaarddiagramtorusmeridional}}    \subcaption{ $\cD(\mu)$}
  \end{minipage}
    \begin{minipage}{0.3\linewidth}\centering
    %\scalebox{1}{\includesvg[pretex=\relscale{1}]{heegaarddiagramtoruslongitudinal}}    \subcaption{ $\cD(\lambda)$}
  \end{minipage}

  \caption{Bordered-sutured Heegaard diagrams for $T^2 \times I$. \label{fig:torusdiagrams}}%

\end{figure}
\begin{prop}The type $D$ structure $\bsd(\cD(\mu))$ is freely generated over $\FF_2$ by a single generator $a$ satisfying $\iota_0 a \iota_0 = a$. Similarly, the type $D$ structure $\bsd(\cD({\lambda}))$ is generated by a single generator $b$ with $\iota_1 b \iota_1 = b$ respectively. In both cases, the map $\delta^1$ is trivial.
\end{prop}
\begin{proof}There are no regions in these diagrams which are not adjacent to some marked point, so the differential is trivial.
\end{proof}

\subsubsection{Morphisms}
\begin{figure}
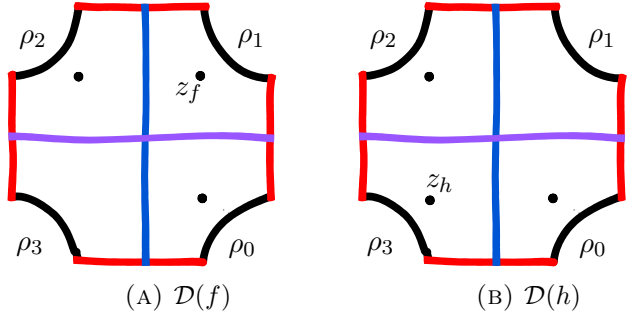

  \begin{minipage}{0.3\linewidth}\centering
    \scalebox{1}{\includesvg[pretex=\relscale{1}]{heegaarddiagramtorustriplef}}    \subcaption{ $\cD(f)$}
  \end{minipage}
    \begin{minipage}{0.3\linewidth}\centering
    \scalebox{1}{\includesvg[pretex=\relscale{1}]{heegaarddiagramtorustripleh}}    \subcaption{ $\cD(h)$}
  \end{minipage}
  \caption{Marked bordered-sutured Heegaard triples for the maps $D_f$ and $D_h$. \label{fig:torusmorphismdiagrams}}%

\end{figure}

We now choose auxiliary data corresponding with the morphisms $f$ and $h$. Here, we have essentially one choice of the pair $\cG(f)$ and $\cG(h)$: we define it to be the pair $(\cD(\mu), \cD(\lambda))$, or, equivalently, the Heegaard triple $\cD(\mu, \lambda)$ shown in Figure \ref{fig:torusmorphismdiagrams}.

Associated to this we must define morphisms of type $D$ structures $$F(f), F(h): \bsd(\cD(\mu)) \rightarrow \bsd(\cD(\lambda)),$$
which for concreteness we relabel by $d_f := F(f), d_h:= F(h)$. We define these to be the triangle counting maps corresponding to the triple $\cD(\mu, \lambda)$ which avoid base-points $z_f$ and $z_h$ respectively, as in Figure \ref{fig:torusmorphismdiagrams}.

\begin{prop}The morphisms $d_f, d_h: \bsd(\cH_\mu) \rightarrow \bsd(\cH_\lambda)$ are given by
\begin{align*}d_f(a) &= \rho_1 \otimes b\\
d_h(a) &= \rho_3 \otimes b.\end{align*}
\end{prop}
\begin{proof}There is only one region (shown in the figures) in each of these diagrams which is not adjacent to a marked point and has correct orientation, and these regions have the corresponding Reeb chords.
\end{proof}

\subsubsection{Carriage returns}
Finally, we have to construct a chain map from $\cone(d_f)$ to $\cone(d_h)$. In this scenario it helps to observe that:
\begin{prop}There are diagrams $\cD_f(\lambda + \mu)$ and $\cD_h(\lambda + \mu)$ such that
 \begin{align*}\cone(d_f) &= \bsd(\cD_f(\lambda + \mu)),\\\cone(d_h) &= \bsd(\cD_h(\lambda + \mu)).\end{align*}
\end{prop}
\begin{proof}
\begin{figure}
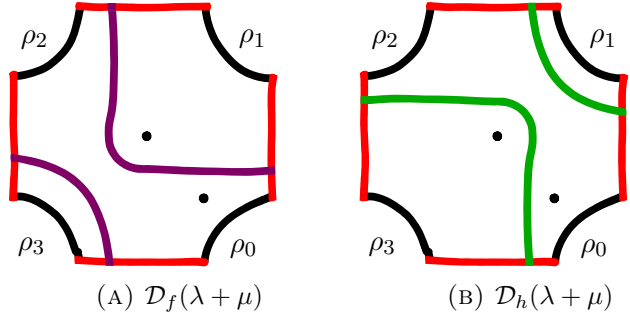

  \begin{minipage}{0.3\linewidth}\centering
    \scalebox{1}{\includesvg[pretex=\relscale{1}]{heegaarddiagramtorusbypassedf}}    \subcaption{ $\cD_f(\lambda + \mu)$}
  \end{minipage}
    \begin{minipage}{0.3\linewidth}\centering
    \scalebox{1}{\includesvg[pretex=\relscale{1}]{heegaarddiagramtorusbypassedh}}    \subcaption{$\cD_h(\lambda + \mu)$}
  \end{minipage}

  \caption{Heegaard diagrams representing the cones of $d_f$ and $d_h$. \label{fig:torusdiagramstwisted}}%
\end{figure}

\label{conesarediagrams}
We choose the diagrams as in Figure \ref{fig:torusdiagramstwisted}. It is straightforward to see that the type $D$ structure $\bsd(\cD_f(\lambda + \mu))$ has two generators $p, q$ with $\delta^1(p) = \rho_1 \otimes q$, so that the map defined by $(a,0) \mapsto p$ and $(0, b) \mapsto q$ is an isomorphism of type $D$ structures. Similarly, $\bsd(\cD_h(\lambda + \mu))$ has two generators $r, s$ with $\delta^1(r) = \rho_3 \otimes s$, and the map $(a, 0) \mapsto r$ and $(0, b) \mapsto s$ is an isomorphism.
\end{proof}
We can therefore work with $\bsd(\cD_f(\lambda + \mu))$ and $\bsd(\cD_h(\lambda + \mu))$ as models for the cones --- this is a setting where we have better means to define maps. Indeed:

\begin{prop}There is a morphism of type $D$ structures $d_{CR}: \cone(d_f) \rightarrow \cone(d_h)$, defined by
$$d_{CR}(a, b) := (\rho_{2} \otimes b, \rho_{123}\otimes a).$$
\end{prop}
\begin{proof}Under the aforementioned isomorphisms from Proposition \ref{conesarediagrams}, the map is defined by counting triangles in the bordered Heegaard triple shown in Figure \ref{fig:toruswitness}.
\begin{figure}
  \begin{minipage}{0.3\linewidth}\centering
    \scalebox{1}{\includesvg[pretex=\relscale{1}]{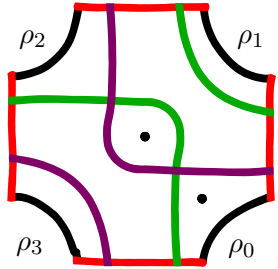}}  
  \end{minipage}
  \caption{A bordered-sutured Heegaard triple $\cD(f,h)$ for the map $d_{CR}$. \label{fig:toruswitness}}%

\end{figure}
There are two triangle shaped domains which support holomorphic representatives; it is straightforward to see that these have the Reeb chords specified.
\end{proof}

WIth these in hand, we have a typewriter in $\mathsf{T}{}^{\cA(T^2)}\mathsf{Mod}$ given by $${\bf m} : = (\bsd(\cD(\mu)), \bsd(\cD(\mu)); d_f, d_h; d_{CR}).$$

\begin{prop}The $DD$ bimodule $\dify({\bf m})$ is isomorphic to $\cfdd(\II)$, as defined in \cite[Figure A.6]{LOT}.
\end{prop}
\begin{proof}Recall that the bimodule $\cfdd(\II)$ is two-dimensional, with generators $x$ and $y$. The structure maps are
\begin{align*}\delta^1(x) &= (\rho_1 \otimes \sigma_3 + \rho_3 \otimes \sigma_1 + \rho_{123} \otimes \sigma_{123}) \otimes y\\
\delta^1(y) & = \rho_2 \otimes \sigma_2 \otimes x.\end{align*}

Similarly, $\dify({\bf m})$ has two generators $a$ and $b$, with structure map
\begin{align*}\delta^1(a) &= (\rho_1 \otimes h + \rho_3  \otimes f + \rho_{123} \otimes fgh) \otimes y\\
\delta^1(b) & = \rho_2 \otimes g \otimes x.\end{align*}
\end{proof}

\subsection{Weak homotopy invariants for bordered-sutured three-manifolds}
Our aim here is to define a set of carriage returns for the weak invariant $\bsd$, ultimately concluding that this gives a weak homotopy invariant of $\diffborsuttor_2^\perp$.

\subsubsection{Auxiliary data}
We start by defining sets of auxiliary data $\cG$ that we will use.

First, for for objects: let $Y$ be an object of $\contdiffborsut$, and take the data $\cG(Y)$ to be the set of all admissible bordered-sutured Heegaard diagrams for $Y$, with edges dictated by requiring an edge between diagrams $\cD, \cD'$ for each handleslide, isotopy or stabilisation relating $\cD$ to $\cD'$.

The second data we need is for morphisms. To any pair of orthogonal objects $Y$ and $Y'$ we define $\cV(f)$ and $\cV(h)$ to be the same subsets of $\cV(Y) \times \cV(Y')$, given by the set of diagrams such that there exists a bordered-sutured Heegaard diagram $$\bcD := \bcD(\cD_1, \cD_2)$$ for the bordered-sutured three-manifold $\sunn(Y, f)$ with
\begin{align*}\cD_1 &= \bcD(\cD_1, \cD_2) \cup \cD(\mu),\\
\cD_2 &=\bcD(\cD_1, \cD_2) \cup \cD(\lambda).\end{align*}

For edges we note that if $\bpsi$ is a handleslide, isotopy or stabilisation move between bordered-sutured diagrams $\bcD$ and $\bcD'$, and $\cD := \bcD \cup \widehat{\cD}$ and $\cD' := \bcD' \cup \widehat{\cD}$ then there is a corresponding handleslide, isotopy or stabilisation $\psi: \cD \rightarrow \cD'$ given by $\psi = \bpsi \cup \mathrm{Id}$, i.e. by extending by the identity on $\widehat{\cD}$.

So, we dictate that if
\begin{align*}(\cD_1, \cD_2) &=  (\bcD(\cD_1, \cD_2) \cup \cD(\mu), \bcD(\cD_1, \cD_2) \cup \cD(\lambda)),\\
(\cD'_1, \cD'_2) &=  (\bcD(\cD'_1, \cD'_2) \cup \cD(\mu), \bcD(\cD'_1, \cD'_2) \cup \cD(\lambda))\end{align*}
then the set of edges in $\cE(f)$ connecting $(\cD_1, \cD_2)$ and $(\cD_1', \cD_2')$ is given by the set of pairs $(p_1, p_2) \subset \cE(Y) \times \cE(Y')$ which are given by extending a handleslide, stabilisation and isotopy from $\bcD(\cD_1, \cD_2)$ to $\bcD(\cD'_1, \cD'_2)$.

Finally, we desire data for carriage returns. This is greatly simplified by coincidence of $\cG(f)$ and $\cG(h)$; we define the data $\cG(f, h)$ to be the diagonal of the product $\cG(f) \times \cG(h)$.

\subsubsection{The invariant}\label{weakinvariantofpairsborderedsutured}
To each of these pieces of auxiliary data we must associate the corresponding constructions in $\amod$. We will need an auxiliary lemma:
\begin{lem}\label{pairinginducesidentifications}Let $(\cD_1, \cD_2) \in \cV(f)$ (resp. $(\cD_1, \cD_2) \in \cV(h)$). Then \begin{align*}\bsd(\cD_1) &= \bsda(\bcD(\cD_1, \cD_2)) \bbox \bsd(\cD(\mu))\\ 
\bsd(\cD_2) &= \bsda(\bcD(\cD_1, \cD_2)) \bbox \bsd(\cD(\lambda)),\end{align*}
where the symbol `$=$' denotes isomorphism.

Furthermore, if $(\cD_1, \cD_2)$ and $(\cD_1', \cD_2')$ are both in $\cV(f)$ (resp. $(\cD_1, \cD_2)$ and $(\cD_1', \cD_2')$ are both in $\cV(h)$) and $\bcD(\cD_1', \cD_2')$ differs from $\bcD(\cD_1, \cD_2)$ by a handleslide, isotopy or stabilisation $\bpsi$ with corresponding invariance map $$\bsda(\bpsi):= \bsda(\cD(\cD_1, \cD_2)) \rightarrow \bsda(\cD(\cD'_1, \cD'_2)),$$ then the corresponding maps
\begin{align*}\psi_1 := \bsda(\bpsi) \bbox \II: \bsda(\cD(\cD_1, \cD_2)) &\bbox \bsd(\cD(\mu)) \rightarrow \bsda(\cD(\cD'_1, \cD'_2)) \bbox \bsd(\cD(\mu)) \\
\psi_2 :=\bsda(\bpsi) \bbox \II: \bsda(\cD(\cD_1, \cD_2)) &\bbox \bsd(\cD(\lambda)) \rightarrow \bsda(\cD(\cD'_1, \cD'_2)) \bbox \bsd(\cD(\lambda)) 
\end{align*}
agree with the maps $\bsd(\bpsi \cup \mathrm{Id(\mu)})$ and $\bsd(\bpsi \cup \mathrm{Id(\lambda)})$ respectively.
\end{lem}
\begin{proof}The pairing theorems of \cite{LOT, zarev:borderedsutured} give the first half of the result up to homotopy equivalence; the stronger result follows from observing that in both cases, the relevant gluing results used are particularly simple as they do not involve any non-provincial holomorphic curves --- and the result therefore holds up to isomorphism.

For the second half, the result is clear by construction for an isotopy or stabilisation; for a handleslide, one a priori needs gluing results of a form not present in current literature but again, the location of the marked points in $\cD(\mu)$ and $\cD(\lambda)$ preclude any non-provincial holomorphic polygons arising --- cf.  \cite[Proposition 5.34]{lot:ssbdcii}.
\end{proof}

To a diagram $\cD \in \cV(Y)$ we associate the bordered-sutured Floer invariant $\bsd(\cD)$, and to any edge $e \in \cE(Y)$ we associate the corresponding invariance map on bordered-sutured Floer invariants; we will address the remaining edge data shortly.

By the first half of Lemma \ref{pairinginducesidentifications}, it makes sense to associate to any pair $(\cD_1, \cD_2) \in \cV(f)$ and $(\cD_1, \cD_2) \in \cV(h)$ the morphisms $\bsd_f(\cD_1, \cD_2) := \II \bbox d_f$ and $\bsd_f(\cD_1, \cD_2) := \II \bbox d_h$ respectively, under the above identification.

Finally, we address the carriage returns and, for any $(\cD_1, \cD_2; \cD_1, \cD_2) \in \cV(f, h)$ we let $\bsd(\cD_1, \cD_2; \cD_1, \cD_2)$ be the morphism from
$$\cone(\II \bbox d_f) \rightarrow \cone(\II \bbox d_h)$$
given by $\II \bbox d_{CR}$.

This gives a typewriter $\bM(\cD_1, \cD_2)$ for every pair $(\cD_1, \cD_2) \in \cV(f) = \cV(h)$.

There is a more straightforward way to sum up the construction in this section:
\begin{prop}Let $\bM = (M_0, M_1; D_f, D_h; D_{CR})$ be a typewriter in ${}^{\cA_2}\mathsf{Mod}$. Then $\cdot \bbox \bM$ constitutes a functor from the category ${}^{\cA_1}\mathsf{Mod}_{\cA_2}$ to the category ${}^{\cA_1}\mathsf{Mod}$ which descends to the homotopy category.

\end{prop}
\begin{proof}This is a straightforward application of functorial properties of the box tensor product $\bbox$.
\end{proof}

The typewriter $\bM(\cD_1, \cD_2)$ is then equal, by construction, to the typewriter $$\bsda(\bcD(\cD_1, \cD_2))\bbox {\bf m}.$$

This also allows us to associate the data for edges which we need: if $(p_1, p_2) \in \cE(f) = \cE(h)$ connects the pair $(\cD_1, \cD_2)$ with $(\cD_1' ,\cD_2')$ then there is a corresponding path $\bp(p_1, p_2)$ connecting $\bcD(\cD_1, \cD_2)$ with $\bcD(\cD_1' ,\cD_2')$; we let the corresponding morphism of typewriters be given by
$$\bsda(\bp(p_1, p_2)) \bbox \II: \bsda(\bcD(\cD_1, \cD_2)) \bbox  {\bf m} \rightarrow \bsda(\bcD(\cD'_1, \cD'_2)) \bbox  {\bf m}$$
inducing the requisite edge maps by Proposition \ref{carriagereturnsistypewriters} (using Lemma \ref{pairinginducesidentifications} to ensure compatibility with the underlying edge maps $\bsd(p_i)$).

\begin{rmk}We will \emph{not} construct carriage advances. There is an inherent problem in doing so with the above set-up: at time of writing, there is no appropriate pairing theorem in the literature for extendable or partially extendable type $D$ structures.
\end{rmk}

In summary, we have constructed a set of carriage returns for the weak invariant $\bsd$, extending it to a  weak homotopy invariant of $\contdiffborsut$ with carriage returns, valued in $\amod$.

As a result:
\begin{thm}There is a weak homotopy invariant of $\diffborsuttor_2^\perp$ which we call $\mathsf{T}\bsd$, valued in the category $\mathsf{T}\amod$. In particular, this restricts to the subcategory $\diffsuttor_2^\perp$ to give a weak homotopy invariant $\mathsf{T}SFH$.
\end{thm}
\begin{proof}This is a direct consequence of Proposition \ref{carriagereturnsistypewriters}.
\end{proof}

\subsubsection{Equivalence}
It turns out that in this instance, we recover the bordered-sutured invariants proper:
\begin{thm}The weak homotopy invariant $\mathsf{T}\bsd$ satisfies $$\dify \mathsf{T}\bsd(Y, f) = \bsdd(\sunn(Y,f)),$$
and $\mathsf{T}SFH$ satisfies $$\dify \mathsf{T}SFH(Y, \gamma) = \bsd(\sunn(Y, f).$$
\end{thm}
\begin{proof}For any pair $(\cD_1, \cD_2)$ in $\cV(f) = \cV(h)$ we have that
\begin{align*}\dify(\bM(\cD_1, \cD_2)) &= \dify(\bsda(\bcD(\cD_1, \cD_2)) \bbox {\bf m})\\
& = \bsda(\bcD(\cD_1, \cD_2)) \bbox \dify({\bf m}) \\
& = \bsdd(\bcD(\cD_1, \cD_2)) \bbox \cfdd(\id)\\
&\cong \bsdd(\bcD(\cD_1, \cD_2)).\end{align*}\end{proof}

\subsubsection{Relationship with sutured constructions}
The above section feels somewhat circular: we define the values of $\bsd(\cD_1, \cD_2)$ using bordered-sutured Floer invariants themselves. It is possible to attempt to circumvent this by defining, for any $(\cD_1, \cD_2)$ as specified, the values of $F_f(\cD_1, \cD_2)$ and $F_h(\cD_1, \cD_2)$ by counting holomorphic triangles in the bordered-sutured Heegaard triple given by $\bcD(\cD_1, \cD_2) \cup \cD(\mu, \lambda)$ --- this gives morphisms which are homotopic to the ones used above, by \cite[Proposition 5.35]{lot:ssbdcii}. In particular, when restricting attention to the purely sutured subcategory this allows one to just count triangles in sutured diagrams.

Unfortunately at this point it ceases to be straightforward how to proceed. Before, we had an explicit identification of the mapping cone of both $\II \bbox d_f$ with $$\bsda(\cD(\cD_1, \cD_2)) \bbox \cD_f(\lambda + \mu)$$ and  $\II \bbox d_h$ with $$\bsda(\bcD(\cD_1, \cD_2)) \bbox \cD_h(\lambda + \mu),$$ allowing us to define a carriage return by counting triangles in a bordered Heegaard triple. Of course, we can still use \cite[Proposition 5.35]{lot:ssbdcii} to identify $\II \bbox d_{CR}$ with a count of triangles in the Heegaard triple $\bcD(\cD_1, \cD_1) \cup \cD_{CR}$ --- but it is only true that, for instance, $\bsda(\bcD(\cD_1, \cD_2)) \bbox \bsd(\cD_f(\lambda + \mu))$ is \emph{homotopy equivalent} to the cone of $\bsd_f(\cD_1, \cD_2)$, so we do not a priori obtain a carriage return from this triangle count.

It's presumably possible to track things through the equivalences induced by pairing, as in \cite[Theorem 8]{lot:ssbdcii}--- but at this point one is probably no better off than using the construction given in the previous subsection.

A more fruitful avenue of approach is to instead construct the data for a triangulated weak homotopy invariant, and show that it admits carriages. This approach will be discussed further in \cite{me:iteratedcones}.
\bibliographystyle{plain}
\bibliography{references}
\end{document}